{
\documentclass[11pt]{article}\usepackage{url}
\usepackage{enumerate}
\usepackage{graphicx} 
\usepackage{float} 
\usepackage{amssymb}
\usepackage{amsmath}
\usepackage{mathtools}
\usepackage{txfonts}
\usepackage{stmaryrd}

\usepackage{color}

\def\rr{{\Bbb R}}
\def\rz{{{\rr}^n}}
\def\zc{{\Bbb Z}}

\def\cc{{\Bbb C}}

\def\ccl{{\Bbb C}^{n+1}}

\def\f{\mathbf{f}}
\def\L{\mathcal{L}}
\def\A{\mathcal{A}}
\def\W{\mathcal{W}}

\def\R{\mathcal{R}}

\def\rdd{\rr_+^{n+1}}
\def\rdl{\rr^{n+1}}
\def\fz{\infty}

\def\supp{{\rm{\ supp\ }}}
\def\loc{{\rm{\ loc\ }}}

\def\ez{\epsilon}

\def\supp{{\rm supp}}
\def\loc{{\rm loc}}

\def\l{\left}
\def\r{\right}

\def\dsum{\displaystyle\sum}

\def\dint{\displaystyle\int}

\def\dsup{\displaystyle\sup}

\newcommand{\rn}{\mathbb{R}^n}

\newcommand{\reu}{\mathbb{R}^{n+1}_+}

\usepackage{theorem}

\newtheorem{theorem}{Theorem}[section]
\newtheorem{lemma}{Lemma}[section]

{
    \theoremheaderfont{\bfseries}
    \theorembodyfont{\normalfont}
}
\graphicspath{{figure/}}                                 
\usepackage[a4paper]{geometry}                           
\geometry{left=2.5cm,right=2.5cm,top=2.5cm,bottom=2.5cm} 
\linespread{1.2}                                         

\begin{document}
\title{\bf  $L^2$ estimates for commutators of the Dirichlet-to-Neumann Map 
associated to elliptic operators with complex-valued bounded measurable coefficients on $\rdd$} 
\date{}
\author{\sffamily Steve Hofmann, Guoming Zhang}
\renewcommand{\thefootnote}{\fnsymbol{footnote}}
\footnotetext[1]{The first author is supported  by NSF grant DMS-2000048. 
The second author is supported  by the NNSF (11771023) of China.}

\parskip=2pt

\maketitle

{\noindent{\bf Abstract:}
     In this paper we establish 
     commmutator estimates for the Dirichlet-to-Neumann Map associated to a divergence form elliptic operator 
     in the upper half-space $\reu:=\{(x,t)\in \rn \times (0,\infty)\}$,
     with uniformly complex elliptic, $L^{\fz}$, $t$-independent 
     coefficients.   By a standard pull-back mechanism,
these results extend corresponding results of Kenig, Lin and Shen for the Laplacian in a Lipschitz domain,
which have application to the theory of homogenization.

\vspace{1ex}

\tableofcontents

\section{Introduction}

Let $\L :=-\mbox{div}(A\nabla )$, defined in $\rdl(n\geq1)$,  where $A=A(x)$ is a $n+1\times n+1$ matrix with complex-valued, bounded and $t$-independent coefficients satisfying the uniform (complex)-ellipticity condition
$$\gamma|\xi|^{2}\leq \mbox{Re}\left \langle A(x)\xi, \xi \right\rangle =\mbox{Re}\l( \sum_{i, j=1}^{n+1}A_{ij}(x)\xi_i\overline{\xi_j}\r),  \quad\quad\left\|A\right\|_{L^{\fz}}\leq \gamma^{-1}, \eqno(1.1)$$
 for some $\gamma \in (0, 1],$ and all $\xi\in \ccl,\; x\in \rz.$  Moreover, throughout our paper, we shall further assume that there exists $A_0$, a uniformly elliptic, $t$-independent matrix as above, which in addition is real and symmetric, such that
$$\|A-A_{0}\|_{L^{\fz}}\leq \ez,\eqno(1.2)$$ where $\ez$ depends only on $n, \gamma$.

If we assume $f \in C^{\fz}_0(\rz)$, then the Dirichlet problem  $$\l\{
           \begin{aligned}
           \L u=0\quad\;\mbox{in}\;\rdd &\\
           \lim_{t\rightarrow 0}u(\cdot,t)=f,\quad\; &
           \end{aligned}\r.\eqno(1.3) $$ has a unique solution $u\in \dot{W}^{1,2}(\rdd)$, the 
           space of functions modulo constants with seminorm given by the norm of $\nabla u$ 
           in $L^2(\rdd)$,
           and the Dirichlet-to-Neumann map,
           defined by 
           $$f\to \Lambda(f):=\frac{\partial{u}}{\partial{\nu_A}}
           =\partial_{\nu_A}u:=-e_{n+1}\cdot A(\nabla u)\biggl|_{t=0}\,,$$ 
           extends to
           a mapping from  $\dot{H}^{1/2}(\rn)\cong I_{1/2}(L^2(\rn))$ to
           $\dot{H}^{-1/2}(\rz)$, where $\dot{H}^{-1/2}(\rz)$ denotes the dual space of the 
           fractional Sobolev space $\dot{H}^{1/2}(\rz)$ (see \cite{HK}); here $I_{1/2}$ denotes the usual 1/2
           order homogeneous fractional integral operator (i.e., Riesz potential). 
           We also define the homogeneous Sobolev space $\dot{L}^2_1(\rz)$ to be the completion of $C^{\fz}_0(\rz)$ with respect to the seminorm $\|\nabla f\|_{2}$.   For convenience, we set $\dot{H}^1(\rn) := \dot{L}^2_1(\rz)$, and we define the
           inhomogeneous version by
      $H^{1}(\rz)=L^{2}(\rz)\cap \dot{H}^1(\rn)$. 
For $\epsilon>0$ small enough, depending only on $n$ and $\gamma$,
we obtain that 
$$\Lambda: \dot{H}^1(\rz) \rightarrow L^2(\rz).\eqno(1.4)\,.$$ 
           In fact, (1.4) is an immediate consequence of the solution of the Regularity problem given in
           \cite[Theorem 1.14]{AAA}. 
           
           We let $C^{0,1}(\rn)$ denote the space of Lipschitz functions, with norm
\[
\|g\|_{C^{0,1}(\rn)}   := \sup_{x,y\in \rn:\,x\neq y } \frac{|g(x)-g(y)|}{|x-y|}  \,.
\]     
           
Now we can state our main results as follows.   The first generalizes the classical commutator theorem of 
A. P. Calder\'on \cite{Ca}.

\begin{theorem}
Suppose that $A$ satisfies (1.1) and (1.2) with $\ez$ sufficiently small,
depending on dimension and ellipticity. 
Then, for any $f\in L^2(\rz)$ and $g\in C^{0, 1}(\rz),$  $$\|[\Lambda, g](f)\|_{L^{2}(\rz)}\leq C\|f\|_{L^{2}(\rz)}\|g\|_{C^{0, 1}(\rz)},\eqno(1.5)$$ where the constant $C$ depends only on $n$ and $\gamma$.
\end{theorem}

\begin{theorem}
Suppose that $A$ satisfies the hypotheses of Theorem 1.1. Then, for any 
$f\in L^{\fz}(\rz)$ and $g\in \dot{H}^{1}(\rz),$  $$\|[\Lambda, g](f)\|_{L^{2}(\rz)}\leq C\|f\|_{L^{\fz}(\rz)}\|g\|_{\dot{H}^{1}(\rz)},\eqno(1.6)$$ where the constant $C$ depends only on $n$ and $\gamma$.
\end{theorem}

Analogous results were previously obtained in \cite{KLS}, for the Laplacian
in a Lipschitz domain, as part of the authors' study of homogenization.
 Our results, with $A$ real and symmetric, include this case, by a well-known pullback
mechanism.  A different
generalization of the results in \cite{KLS} has been obtained in  \cite{SZ}, for elliptic systems in
Lipschitz domains, with H\"older continuous coefficients.  Neither our work nor that of \cite{SZ} subsumes the
other.  The approach to these commutator results in both papers \cite{KLS} and \cite{SZ}, is based on
a bilinear estimate of Dahlberg \cite{D}, and it's extension in \cite{SZ} to certain variable coefficient elliptic systems.
In \cite{HS1}, the first author of this paper established Dahlberg's bilinear estimate 
for the class of second-order elliptic operators enjoying the same assumptions  
that we impose here, i.e., the matrix $A$ satisfies (1.1) and (1.2) with $\ez$ small enough, in
 the upper half-space. The latter result, along with layer potential technology for
 the operators under consideration, will allow us to follow the strategy 
 in \cite{KLS} and \cite{SZ}, to obtain the stated theorems.



\smallskip

\noindent {\bf Remark 1.} In our Theorem 1.2, as compared to its analogue
\cite[Theorem 1.2]{SZ}, we obtain an estimate in terms of
the $L^\infty$ norm of the boundary data, as opposed to that of the solution $u$ itself, since we establish an 
Agmon-Miranda maximum principle for our solutions
(see \cite[Remark 1.4]{SZ}, and Section 4 below).

\smallskip
	
The paper is organized as follows.  In the next section we discuss certain preliminaries.  In Section 3, we prove
Theorem 1.1.  In Section 4, we establish an Agmon-Miranda maximum principle for the class of operators under
consideration, which we then use in Section 5 to give the proof of Theorem 1.2.

\section{Preliminaries}

We begin by setting some notational conventions.  For convenience, we 
often write $B\lesssim D$, and $B\approx D$, to mean that there exists a positive 
constant $C$, depending only on dimension and  the quantitative hypotheses of our theorems, 
such that respectively, $B\leq C D$, and $C^{-1} D \leq B\leq C D$. 
We normally use $Q$ to
denote cubes  in $\rz$, and for $\lambda>0$, we let $\lambda Q$ be the concentric dilate of $Q$
with side length $\lambda \ell(Q)$.

Let us now 
recall the De Giorgi-Nash-Moser estimates:  
under the same assumptions as in Theorem 1.1 (in fact $t$-independence is not required), 
there is a constant $C$ and an exponent $\alpha>0$, both depending only on  dimension and ellipticity, such that
for any ball $B=B(X, R)$, 
if $\L u=0$ in $2B=B(X, 2R)$, 
then $$|u(Y)-u(Z)|\leq C\l(\frac{|Y-Z|}{R}\r)^{\alpha}\l(\fint_{2B}|u|^2\r)^{1/2},
\eqno(\text{DGN})$$ 
whenever $Y, Z\in B$, and 
$$\dsup_{Y\in B}|u(Y)|\leq C\l(\fint_{2B}|u|^2\r)^{1/2}\,.
\eqno(\text{M})$$ 
As is well known, these results
may be found in \cite{D, M, N} in the case of real coefficients; the extension to the case of complex perturbations
of real coefficients is due to Auscher \cite{A} (see Theorem 4.1 below). 

We shall make use of the theory of layer potential operators associated to an operator $\L=-\text{div}\,A \nabla$
as in (1.1), 1.2).  
Let $E(x,t,y,s)=E_\L(x,t,y,s)$ be the fundamental solution of $\L$. The existence 
of the fundamental solution in our setting is given in \cite{HK2}. By 
$t$-independence of our coefficients, 
we have that $$E(x,t,y,s)=E(x,t-s,y,0).\eqno(2.1)$$ The single and double layer 
potential operators associated to $\L$ are defined, respectively, by 
$$\mathcal{S}_tf(x)=\mathcal{S}^\L_tf(x):= \int_{\rz}E(x,t,y,0) f(y) dy, \quad t\in \rr, \quad x\in \rz,$$
$$\mathcal{D}_{t}f(x)=\mathcal{D}^\L_{t}f(x):= \int_{\rz}\overline{\partial_{\nu_{A^*}}E^*(y,0,x,t)} f(y) dy, 
\quad t\neq 0, \quad x\in \rz,\eqno(2.2)$$ where $A^*$ is the hermitian adjoint of $A$ and 
$$\partial_{\nu_{A^*}}E^*(y,0,x,t)
=-\sum_{j=1}^{n+1}A^*_{n+1, j}(y)\frac{\partial{E^*}}{\partial{y_{j}}}(y,0,x,t)=-e_{n+1}\cdot A^*(y)\nabla_{y,s}
E^*(y,s,x,t)\biggl|_{s=0}.\eqno(2.3)$$ 
Here, $E^*=E_{\L^*}$ 
denotes the fundamental solution of $\L^*$, the 
hermitian adjoint of $\L$, and we have $$\overline{E^*(y,s,x,t)}=E(x,t,y,s).\eqno(2.4)$$ We shall 
use the following notations: 
$D_{j}:= \frac{\partial}{\partial{x_j}}=\partial_{x_j},  1\leq j\leq n+1$, where $x_{n+1}:= t$ (so that
$D_{n+1}=\partial_t$), and for
a vector ${\bf v}=(v_1,v_2,...v_{n+1})\in \mathbb{R}^{n+1}$, we let
${\bf v}_\|:= (v_1,v_2,...v_{n},0)\cong (v_1,v_2,...v_{n})$
denote the projection of ${\bf v}$ onto $\rn\times \{0\}$.  
Similarly, we define
$\nabla_{\parallel}:=(\partial_{x_1},..., \partial_{x_n}).$ 
We shall set
$$(\mathcal{S}_t\nabla)\,f(x):= \int_{\rz}\nabla_{y,s}E(x,t,y,s)\bigl|_{s=0}\,f(y)dy,$$ 
so that 
$$
(\mathcal{S}_tD_{n+1})
=-\partial_t \mathcal{S}_t\,, \qquad
(\mathcal{S}_t\nabla_{\parallel}) \,\vec{f}=
-\mathcal{S}_t\,\big(\mbox{div}_{\parallel}\vec{f}\,\,\big)\,,
\eqno(2.5)$$ 
for, say, $\vec{f}\in C^1_0(\rn,\mathbb{C}^n)$.
For all $m\geq 1$, it follows from (2.1)-(2.5) that 
$$\mbox{adj}\big(\nabla \partial^{m-1}_t (\mathcal{S}_{t}\nabla_{\parallel})\big)
=\pm \nabla_{\parallel}\partial^{m-1}_t(\mathcal{S}^{\L^*}_{-t}\nabla)\eqno(2.6)$$ 
(here the choice of ``plus" or ``minus" depends on $m$),
and 
$$\mbox{adj}(\nabla \partial^{m-1}_t \mathcal{D}_{t})
=\pm \partial_{\nu_{A^*}}\partial^{m-1}_t(\mathcal{S}^{\L^*}_{-t}\nabla),\eqno(2.7)$$ 
where $\mathcal{S}^{\L^*}_{t}$ denotes the single layer potential with associated to $\L^*$,
the adjoint co-normal derivative is defined by
$\partial_{\nu_{A^*}}=-\sum_{j=1}^{n+1}A^*_{n+1, j} \,D_j$, and $\mbox{adj}(T)$ denotes the hermitian adjoint of an operator $T$ acting in $\rz$.

\noindent{\bf Remark 2.} By \cite[Theorem]{AAA}, for $\L=\text{div}\,A\nabla$, with $A$ satisfying 
(1.1) and (1.2) with $\epsilon$ small enough depending on dimension and ellipticity,
we have the layer potential bounds\footnote{These bounds continue to hold in the absence of
condition (1.2) (for a suitable definition of the layer potentials):  see \cite{R}.}
\[
\sup_{t\neq 0}\left(\|\nabla\,\mathcal{S}_t^\L f\|_{L^2(\rn)} \, +\, \|(\mathcal{S}_t^\L\,\nabla) f\|_{L^2(\rn)}\right) \, \lesssim \, \| f\|_{L^2(\rn)}\,.
\]
In particular, this yields $L^2$ boundedness of the double layer potential $\mathcal{D}_t$, uniformly in $t$.
Of course, analogous results hold with $\L$ replaced by its adjoint $\L^*$.

Given $x_0\in \rz$ and $\beta>0$, define the cone $\Gamma_{\beta}(x_0):=\{(x,t)\in \rdd; |x_0-x|<\beta t\},$ then for measurable function $F: \rdd\rightarrow \cc$, the non-tangential maximal operator 
$N^{\beta}_*$ is defined $$N^{\beta}_*(F)(x_0):=\sup_{(x,t)\in \Gamma_{\beta}(x_0)}|F(x,t)|\,,$$ and 
note that when $\beta = 1$, we shall often simply write $\Gamma = \Gamma_1$, and
$N_*(F):=N^{1}_*(F)$.
We recall that the $L^2$-norms of $N_*=N^{1}_*$ and $N^{\beta}_*$ are equivalent 
for any $\beta>0$ (see \cite{FS}). 
Following \cite{KP}, we also introduce $$\tilde{N}_*(F)(x_0):=\sup_{(x,t)\in \Gamma_{1}(x_0)}\l(\fint_{|(x,t)-(y,s)|<t/4}|F(y,s)|^2dyds\r)^{1/2},$$ where the symbol $\fint$ denotes the mean value, i.e. $\fint_{E}f\equiv |E|^{-1}\int_{E}f.$\;
We say $u\rightarrow f\; n.t.$ to mean that for a.e. $x\in \rz$, $\lim_{(y,t)\rightarrow (x,0)}u(y,t)=f(x)$, where the limit runs over $(y,t)\in \Gamma(x)$, and in the sequel, we shall use the notation 
$\interleave  \cdot  \interleave $ as a short-hand for the $T_2^2$ tent-space norm (see \cite{CMS}),
i.e.,
$$\interleave  F  \interleave:= \l(\iint_{\rdd}|F(x,t)|^2\frac{dxdt}{t}\r)^{1/2}.$$

Next, we state a technical lemma concerning
 the single layer potential, as well as general solutions.  The lemma will
follow essentially immediately from known results, and will be useful in the sequel.

\begin{lemma}
Suppose that $A,\L$ satisfy the same hypotheses as in Theorem 1.1.
Let $\f\in L^2(\rz,\ccl)$, $f \in L^{\infty}(\rz,\ccl)$, and let $m\geq 1$.  Then
 $$\sup_{t > 0}\|t^m\nabla\partial^{m-1}_t(\mathcal{S}^{\L}_{t}\,\nabla)\, \f\|_{L^2(\rz)}\lesssim
  \l\|\f\r\|_{L^2(\rz)},\eqno(2.8)$$
$$\interleave  t^m\nabla\partial^{m-1}_t(\mathcal{S}^{\L}_{t}\,\nabla)\, \f  \interleave 
\lesssim \|\f\|_{L^2(\rz)}, \eqno(2.9)$$ 
and
$$\sup_{Q}\frac{1}{|Q|}\int^{\ell(Q)}_{0}\int_{Q}\big|t^m\nabla\partial^{m-1}_t
(\mathcal{S}^{\L}_{t}\,\nabla)\, 
\f(x)\big|^2\frac{dxdt}{t}\lesssim \|\f\|^2_{L^{\infty}(\rn)}\,.\eqno(2.10)$$
Furthermore, 
for every cube $Q$ and all $0<t\leq 16\, \ell(Q)$, 
$$\|t^m\nabla\partial^{m-1}_t
(\mathcal{S}^{\L}_{t}\,\nabla)\, (\f1_{2^{k+1}Q\setminus 2^kQ})\|^2_{L^2(Q)}
\lesssim 
2^{-nk}2^{2k}\l(\frac{t}{2^k \ell(Q)}\r)^{2m}\|\f\|^2_{L^2(2^{k+1}Q\setminus 2^kQ)}, \quad 
\forall \,k\geq 1.\eqno(2.11)$$
Finally, suppose that $Lu=0$ in $\mathbb{R}^{n+1}$,
with $\sup_{t>0}\|u(\cdot,t)\|_2 < \infty.$  Then 
$$\||t^m\nabla \partial_t^{m-1} u\|| \lesssim \sup_{t>0}\|u(\cdot,t)\|_2. \eqno(2.12)$$
In all of these estimates, the implicit constants depend on $m$, $n$, and ellipticity.   Of course,
the corresponding estimates hold 
also in the lower half-space, and with $\L$ replaced by $\L^*$.
\end{lemma}

\noindent
{\it Sketch of Proof.}\quad  
Given the $L^2$ bounds discussed in Remark 2,
the case $m=1$ of estimate (2.8) is 
\cite[Lemma 2.11]{AAA}  (we caution the reader that
the exponents in \cite[Lemma 2.11]{AAA} are written differently, so that the case $m=-1$ there 
corresponds to our case $m=1$).  The case $m>1$ may be reduced to the case $m=1$ by an 
induction argument which exploits the ``Caccioppoli on slices" estimate in \cite[Proposition 2.1]{AAA}.
We omit the details.

For $m=1$, the square function bound (2.9) is \cite[Lemma 3.1]{HMM}, the 
Carleson measure estimate (2.10) is \cite[Corollary 3.3, estimate (3.4)]{HMM},
and the square function bound (2.12) is
\cite[Lemma 3.1]{HS1}.  For each of 
(2.9), (2.10), and (2.12), the case 
$m>1$ may  be reduced to the case $m=1$ by an induction argument that uses Caccioppoli's 
inequality in Whitney boxes.  We omit the details.

 Finally, using again the ``Caccioppoli on slices" estimate in \cite[Proposition 2.1]{AAA}, one may
reduce estimate (2.11) to \cite[Lemma 2.9 (i)]{AAA}.  
 Again we omit the details.
\hfill$\Box$

We shall also require some of the main results in \cite{AAA}, which 
we summarize as follows: 

\noindent{{\bf Theorem 2.2}\hspace{0.1cm}{\rm\bf (\cite[Theorem 1.14]{AAA})}. 
 {\em
Suppose that $\L :=-\mbox{div}(A\nabla )$, $A$ and $A_0$ are defined as above, then the Dirichlet problem $$\left\{
\begin{aligned}
\L u=0\;\quad \mbox{in}\; \rdd\quad\quad\quad\quad\quad\quad&\\
\lim_{t\rightarrow 0}u(\cdot,t)=f \quad\mbox{in}\; L^2(\rz)\; \mbox{and}\; n.t.&\\
||N_*(u)||_{2}+ \interleave  t\nabla u  \interleave \lesssim \|f\|_{2},\quad\quad\;\; &
\end{aligned}\right. \eqno(D)^{\L}_2$$ and the Regularity problem $$\left\{\begin{aligned}
            \L u=0\;\quad \mbox{in}\; \rdd\quad\quad\quad\quad\quad\quad&\\
            \lim_{t\rightarrow 0}u(\cdot,t)=f \quad\mbox{in}\; \dot{L}^2_1(\rz)\; n.t.\quad\;\;&\\
            ||\tilde{N}_*(\nabla u)||_{2}\lesssim \|\nabla_{\parallel}f\|_{2},\quad \quad\quad\quad\quad& 
\end{aligned}\right. \eqno(R)^{\L}_2$$ are both solvable if $\ez$ is sufficiently small, depending only on $n$ and $\gamma$. The solution of $(D)^{\L}_2$ is unique and the solutions of $(R)^{\L}_2$ are unique modulo constants.
Analogous conclusions hold for $\L^*$. }

\section{Proof of Theorem 1.1} 

Under the same assumptions of our main theorems, using the results of \cite{AAA}, 
we see that if $u$ is the solution of the Regularity problem $(R)^{\L}_2$, with data
 $f\in \dot{H}^1(\rn)$, then
\begin{equation*}
\|\Lambda f\|_2:=\big\|\frac{\partial{u}}{\partial{\nu_A}}\big\|_2 \lesssim 
\|\tilde{N}_*(\nabla u)\|_2\lesssim \|\nabla _{\parallel}f\|_{2},
\end{equation*}
i.e., (1.4) holds for the Dirichlet-to-Neumann map $\Lambda$.

\noindent{\bf Remark 3.}
We note that for $f\in H^1(\rn)=L^2(\rn)\cap\dot{H}^1(\rn)$, we may solve both
$(R)^{\L}_2$ and $(D)^{\L}_2$ with boundary data $f$, and the respective
resulting solutions $u_R$ amd $u_D$
are the same\footnote{More precisely, they are equal modulo an additive constant.}.
This fact follows from the
``compatible solvability" of the solutions constructed in \cite[Theorem 1.11]{HK}.

The commutator of $\Lambda$ with a function $g$ is defined by 
$$[\Lambda, g](f):=\Lambda(gf)-g\Lambda(f).$$ 
Note that for $gf \in H^1(\rz)$ and $f \in H^1(\rz)$, 
both $\Lambda(gf)$ and $g\Lambda(f) $ are well-defined.  
Let $\varphi\in C_0^{\infty}(B(0, 1))$, 
such that $\varphi$ is radial, $0\leq \varphi\leq 1$ and $\int_{\rz}\varphi=1$, and set 
\[
V(x,t):=P_t \, g:= \varphi_{t}*g\,, \eqno(3.1)
\]
where $\varphi_{t}(x)=t^{-n}\varphi(\frac{x}{t})$. 
Observe that $P_t$ thus defines a 
nice approximate identity.  In particular, $P_t1=1$ and $\nabla P_t 1=0$.
Let us note for future reference the elementary fact that
\[
|\nabla V(x,t)|\,=\,  | \varphi_t * \big(\nabla_\| g\big)(x) \, +\, 
(\partial_t \varphi_t) * (g-C_{x,t}) (x)| \,\lesssim\, 
\fint_{|x-y|<t}|\nabla_\|g(y)|\, dy\,,
\eqno(3.2)
\]
by Poincare's inequality, where we have chosen $C_{x,t}= \fint_{|x-y|<t} g(y)\, dy$,.

For any $h\in C_0^{\infty}(\rz)$, by Theorem 2.2, we let $u$ be the solution of $(D)^{\L}_2$ with boundary data $f$, $u_{fg}$ be the solution of $(D)^{\L}_2$ with boundary data $fg$ and $H$ be the solution of $(D)^{\L^*}_2$ with boundary data $h$. Thus, according to the definition of the Dirichlet-to-Neumann map $\Lambda$, along with a standard variational formulation
of the divergence theorem for solutions,
 $$\begin{array}{cl}
&\displaystyle{\int_{\rz}[\Lambda, g]f\overline{h}=\int_{\rz}\frac{\partial{u_{fg}}}{\partial{\nu}}
\overline{h}-\int_{\rz}\frac{\partial{u}}{\partial{\nu}} g\overline{h}}
\quad\quad\quad\quad\quad\quad\quad\quad\vspace{1.8mm}\\
&\quad\quad\quad\quad\quad\displaystyle{=\iint_{\rdd}A\nabla u_{fg}
 \overline{\nabla H}-\iint_{\rdd}A\nabla u \overline{\nabla H}V-\iint_{\rdd}A\nabla u 
 \nabla V\overline{H}}\vspace{1.8mm}\\
&\quad\quad\quad\quad\quad\displaystyle{=\iint_{\rdd}\nabla u_{fg} 
\overline{A^*\nabla H}-\iint_{\rdd}\nabla u \overline{A^*\nabla H}V-\iint_{\rdd}A\nabla u 
\nabla V\overline{H}}\vspace{1.8mm}\\
&\quad\quad\quad\quad\quad\displaystyle{=
 \iint_{\rdd}\nabla u_{fg} 
\overline{A^*\nabla H}-\iint_{\rdd}\nabla (uV) \overline{A^*\nabla H}
+\iint_{\rdd} u \nabla V \overline{A^*\nabla H}
-\iint_{\rdd}A\nabla u 
\nabla V\overline{H}}\vspace{1.8mm}\\
&\quad\quad\quad\quad\quad\displaystyle{=
\int_{\rz}fg \overline{\frac{\partial{H}}{\partial{\nu_{A^*}}}}
- \int_{\rz}fg \overline{\frac{\partial{H}}{\partial{\nu_{A^*}}}}
+ \iint_{\rdd}u\nabla V\overline{A^*\nabla H}-
\iint_{\rdd}A\nabla u \nabla V\overline{H}}
\vspace{1.8mm}\\
&\quad\quad\quad\quad\quad\displaystyle{=\iint_{\rdd}u\nabla V\overline{A^*\nabla H}-
\iint_{\rdd}A\nabla u \nabla V\overline{H}},
\end{array} \eqno(3.3)$$ 
where in the next-to-last step we have used the fact that 
$u_{fg}(\cdot,0) =fg = u(\cdot,0) V(\cdot,0)$.


We may assume that 
$f\in C_0^{\infty}(\rz)$, by density of the latter space in $L^2(\rn)$. 
Since $g\in C^{0, 1}(\rz)$, if $V$ is defined as above, then $\|\nabla V\|_{L^\infty}\lesssim \|\nabla_\| g\|_{L^\infty}$ and $$d\mu=\vert t\nabla^2 V(x,t)\vert^2\frac{dxdt}{t}\;\mbox{is a Carleson measure on}\; \rdd \;\mbox{with norm}\; \|\mu\|_{\mathit{c}}\lesssim \|\nabla_\| g\|_{L^\infty}.\eqno(3.4)$$ 
Our goal is to show that
$$\bigg\vert\int_{\rz}[\Lambda, g]f \overline{h}\bigg\vert \, \lesssim\,  \|\nabla_\| g\|_{L^\infty}\;
 \|f\|_{2}\; \|h\|_{2}, \quad \forall\; h\in C_0^{\infty}(\rz).\eqno(3.5)$$ 
 To prove (3.5), we see from the equality $(3.3)$ that 
 $$\int_{\rz}[\Lambda, g]f \overline{h}=\iint_{\rdd}u\nabla V\overline{A^*\nabla H}
 \,-\,\iint_{\rdd}A\nabla u \nabla V\overline{H}=:I+J.$$ 
 We observe that the two terms, $I$ and $J$,  are of essentially 
 the same type, since $u$ is the solution of $(D)^{\L}_2$ with boundary data $f\in L^2(\rn)$, 
 while $H$ is the solution of $(D)^{\L^*}_2$ with boundary data $h\in L^2(\rn)$.
We now claim that it suffices to prove the estimate
$$\bigg\vert\iint_{\rdd} \tilde{A}\,\nabla U\cdot  
\W\bigg\vert\lesssim \l(\|N_*(U)\|_{2}+\interleave  
t\nabla U  \interleave\r)\l(\|N_*(\W)\|_{2}+\interleave  t\nabla 
\W \interleave\r),\eqno(3.6)$$ 
where $\tilde{A}$ satisfies (1.1) and (1.2), and where 
$\tilde{\L}\,U:= - \text{div}\,\tilde{A}\nabla\, U = 0$ in $\mathbb{R}^{n+1}_+$.
Indeed, taking (3.6) for granted momentarily,
we may apply the latter estimate to term $J$ in the case that  
the values of $U, \tilde{A}, \W$ are respectively $u, A, \nabla 
 V\overline{H}$, or to term $I$ with these values given respectively 
 by $\overline{H}, \overline{A^*}, \nabla V u$. 
In the former scenario, by Theorem 2.2, we have 
$$\|N_*(u)\|_{2}+\interleave  t\nabla u \interleave \lesssim  \|f\|_{2},$$ 
$$\|N_*(\W)\|_{2}=\|N_*(\nabla V H)\|_{2}\lesssim
 \|\nabla V\|_{L^\infty}\|N_*(H)\|_{2}\lesssim \|\nabla_\| g\|_{L^\infty} \|h\|_{2}\,,$$ 
 and $$\interleave  t\nabla \W \interleave\lesssim \l(\interleave  t\nabla V \nabla H 
  \interleave+\interleave  t\nabla^2 V H\r)\lesssim\l(\|\nabla V\|_{L^\infty}\interleave 
   t\nabla H \interleave+\|\nabla_\| g\|_{L^\infty}\|N_*(H)\|_{2}\r)\lesssim 
   \|\nabla_\| g\|_{L^\infty} \|h\|_{2},$$ where we used (3.4). 
   A similar discussion is applicable to the other case, and (3.5) follows.

It remains to prove (3.6). We will actually prove
a slightly sharper version of (3.6), which is a generalization of the
Dahlberg-type bilinear estimate in \cite{HS1}.  For notational convenience, we shall remove the
``tilde", and just write $\L= -\text{div}\, A\nabla$, where $A$ satisfies (1.1) and (1.2),
and we shall replace $U$ by $u$, and $\W$ by its complex conjugate.
Recall the definition of the standard $Y^{1,2}$ space:
$$Y^{1,2}(\mathbb{R}_+^{n+1}) : = 
\Big\{ u \in  L^{\frac{2(n+1)}{(n+1)-2}}(\mathbb{R}_+^{n+1}): \nabla u \in L^2(\mathbb{R}_+^{n+1})\Big\}.$$

\begin{lemma}
Let $A,\L$ be as above, and let $M$ be an arbitrary bounded
$(n+1)\times(n+1)$ matrix-valued function on $\rn$.
Suppose that $\W\in W^{1,2}_{loc}(\rz, \ccl)$,  and that $u
\in Y^{1,2}(\mathbb{R}_+^{n+1})$ is a solution of
$\L\,u=0 \;\mbox{in} \;\rdd$. 
Then $$\bigg\vert\iint_{\rdd}M\nabla u\cdot  \overline{\W}\bigg\vert\lesssim \,\|M\|_{L^\infty(\rn)}\,
\sup_{t>0}\|u(\cdot, t)\|_{L^2(\rz)}\l(\|N_*(\W)\|_{L^2(\rz)}+
\interleave  t\nabla \W \interleave\r).$$
\end{lemma}
Of course,  Lemma 3.1 (with $u=U$) implies (3.6), since trivially 
$\sup_{t>0}\|u(\cdot, t)\|_{L^2(\rz)} \leq \|N_*(u)\|_{L^2(\rz)}$.

\noindent
{\it Proof.}  Without loss of generality,
we may suppose that $\|M\|_{L^\infty(\rn)}\leq 1$.
The special case that $M=I$, the $(n+1)\times (n+1)$ identity matrix, is proved in
\cite{HS1}, and the argument there may be readily adapted,
{\em mutatis mutandis}, to prove this version.  
For the sake of self-containment, 
and because will shall need to pursue this point anyway to prove Theorem 1.2, we shall
present a slightly different proof (similar but with a small modification) to that of \cite{HS1}.
The proof of Theorem 1.1 will still be a rather routine adaptation of the arguments in
\cite{HS1}.  In the case of Theorem 1.2, matters will be a bit more subtle.
For now, as in \cite{HS1}, it is enough to show
  \begin{multline}\tag{3.7}
 \sup_{0<\rho\ll1}{\bf M}_\rho:=\sup_{0<\rho\ll1}\fint_{\rho}^{2\rho}\bigg\vert\int^{1/\theta}_{\theta}\int_{\rz}
 M(x)\nabla 
 u(x,t)\cdot  \overline{\W(x,t)}dxdt\bigg\vert d\theta \\[4pt]
 \lesssim\,\sup_{t>0}\|u(\cdot, t)\|_{L^2(\rz)}\l(\|N_*(\W)\|_{2}+\interleave  t\nabla \W \interleave\r).
 \end{multline}
We may assume $\sup_{t>0}\|u(\cdot, t)\|_{L^2(\rz)}
+\|N_*(\W)\|_{2}+\interleave  t\nabla \W \interleave<\infty$, since otherwise (3.7) is trivial.  

For each fixed $\rho$,
 integrating by parts in $t$, we obtain the bound
  \begin{multline}\tag{3.8}
  {\bf M}_\rho\leq 
  \fint_{\rho}^{2\rho}\bigg\vert\int^{1/\theta}_{\theta}\int_{\rz}M(x)\nabla 
 \partial_{t}u(x,t)\cdot \overline{\W(x,t)}\,dx\, tdt\bigg\vert d\theta \\[4pt]
 +\quad \fint_{\rho}^{2\rho}\bigg\vert\int^{1/\theta}_{\theta}\int_{\rz}M(x)\nabla u(x,t)\cdot 
\overline{\partial_{t}\W(x,t)}\,dx\,tdt\bigg\vert d\theta
\,\,+\, \, \mbox{``B"} \\[4pt]
=: \, I \,+\, II\,+\, \mbox{``B"} . 
\end{multline}
Here, the term ``B" is 
a sum of two boundary terms, 
satisfying
$$\text{``B"}\,\lesssim\, 
\l(\sup_{r>0}\fint^{2r}_{r}\int_{\rz}r^2
\vert\nabla u(x,t)\vert^2dxdt\r)^{1/2}\sup_{r>0}\|\W(\cdot, r)\|_{2}\lesssim 
\sup_{r>0}\|u(\cdot, r)\|_{2}\sup_{r>0}\|\W(\cdot, r)\|_{2},$$ 
as desired,
where in the integral involving $\nabla u$, we have split $\rz$ into 
cubes of side length $\approx\; r$ and used Caccioppoli's inequality. 
For term $II$ in (3.8), by Cauchy-Schwarz we have the bound
\[ II\,\lesssim\,  \interleave  t\nabla u \interleave
\interleave  t\nabla \W \interleave \, \lesssim \,\sup_{t>0}\|u(\cdot, t)\|_{L^2(\rz)}
 \interleave  t\nabla \W \interleave,
\]
where we have used the case $m=1$ of (2.12) to obtain the inequality 
$ \interleave  t\nabla u \interleave\lesssim \sup_{t>0}\|u(\cdot, t)\|_{L^2(\rz)}$.

We turn now to term $I$.  By
the change of variable $t \rightarrow 2t$, and an integration by parts in $t$, 
\begin{multline*}
I\, =\, 2
\fint_{\rho}^{2\rho}\bigg\vert\int^{1/(2\theta)}_{\theta/2}\int_{\rz}
M(x)\nabla \partial_{t}u(x,2t)\cdot \overline{\W(x,2t)}\,dx\,tdt\bigg\vert d\theta\\[4pt]
\lesssim \, \, 
\fint_{\rho}^{2\rho}\bigg\vert\int^{1/(2\theta)}_{\theta/2}\int_{\rz}
M\nabla \partial^2_{t}u(\cdot,2t)\cdot \overline{\W(\cdot,2t)}\,dx\,t^2dt\bigg\vert d\theta
\\[4pt]
+\,\, \fint_{\rho}^{2\rho}\bigg\vert\int^{1/(2\theta)}_{\theta/2}\int_{\rz}
M\nabla \partial_{t}u(\cdot,2t)\cdot \overline{\partial_{t}\W(\cdot,2t)}\,dx\,t^2dt\bigg\vert d\theta
\,=:\, I_1 + I_2\,.
\end{multline*}
We handle $I_2$ like term $II$ above: by Cauchy-Schwarz, and the case $m=2$ of (2.12), we have
\[
I_2 \lesssim  \interleave  t^2\nabla \partial_t u \interleave  \interleave  t\nabla \W \interleave
\, \lesssim \,\sup_{t>0}\|u(\cdot, t)\|_{L^2(\rz)}
 \interleave  t\nabla \W \interleave\,.
\]

To bound term $I_1$, we integrate by parts again\footnote{The point is to accumulate 
enough $t$-derivatives in order to ensure sufficient decay; see (3.12) below.} 
in $t$, to obtain
\[
I_1\,\lesssim \,I_1' \,+\,I_1''\,,
\]
where
$$I_1'\, =\, 
\fint_{\rho}^{2\rho}\bigg\vert\int^{1/(2\theta)}_{\theta/2}\int_{\rz}
M(x)\nabla \partial^3_{t}u(x,2t)\cdot \overline{\W(x,2t)}\,dx\,t^3dt\bigg\vert d\theta\,,\eqno(3.9)$$
and
\begin{equation*}
I_1'' \, =\,
 \fint_{\rho}^{2\rho}\bigg\vert\int^{1/(2\theta)}_{\theta/2}\int_{\rz}
M\nabla \partial^2_{t}u(\cdot,2t)\cdot \overline{\partial_{t}\W(\cdot,2t)}\,dx\,t^3dt\bigg\vert d\theta \,.
\end{equation*}
Term $I_1''$ can be treated just like terms $II$ and $I_2$, using Cauchy-Schwarz
and the case $m=3$ of (2.12).  We omit the now familiar details.

It remains now only to consider $I_1'$.
 With $t>0$ fixed, set 
 $f_t(x)=\partial_{t}u(x,t)$, and note that 
 \[
 \|f_t\|_{L^2(\rn)} \lesssim t^{-1}\sup_{\tau>0}\|u(\cdot, \tau)\|_{L^2(\rz)}<\infty\,,
 \]
by the Moser-type local boundedness assumption (M) above, and Cacioppoli's inequality.
Moreover, by $t$-independence of the coefficients, $\partial_tu(\cdot,t+\cdot)$ is a solution
in $\rdd$, so by Green's formula\footnote{See, e.g., \cite[Theorem 4.16]{BHLMP} for a justification of
the Green formula in this setting (in fact, in a more general setting).},  
we can write
$$\partial_{t}u(\cdot, t+s)
\,=\,-\mathcal{D}_s(f_t)\,+\,\mathcal{S}_s(\partial_{\nu_A}(\partial_{t}u(\cdot, t+\cdot)).
\eqno(3.10)$$ 
Observe that, at least formally, using $t$-independence and the fact 
that $\partial_tu$ is a solution,
\begin{multline*}
\partial_{\nu_A}\l(\partial_{t}u(\cdot, t+\cdot)\r)
\,=\,-\sum^{n+1}_{j=1}D_{n+1}\bigg(A_{n+1, j}\,D_j
u(\cdot,t+s)\bigg)\big|_{s=0} \\[4pt]
=\, \sum^{n}_{i=1} \sum^{n+1}_{j=1}
D_i\bigg(A_{i,j}D_ju(x, t+s)\bigg)\big|_{s=0}
\,= \,\nabla_{\parallel}\cdot \big(A\nabla u(\cdot, t)\big)_{\parallel}\,,
\end{multline*}
where we interpret the identity in the weak sense on $\rn$,
 see \cite[Lemma 2.15]{AAA}.
Consequently,
setting $s=t$ in (3.10), we get 
\[
\left(\nabla \partial^3_{t}u\right)(\cdot,2t)=-\left(\nabla  \partial^2_{t}\mathcal{D}_t\right)(f_t)\,-\,
\left(\nabla  \partial^2_{t} \big(\mathcal{S}_t\nabla_{\parallel}\big)\right) 
\left(\big(A\nabla u(\cdot, t)\big)_{\parallel}\right)\,, \eqno(3.11)
\]
where we have used (2.5).    We may then obtain the bound
\[
 \left|\int_{\rz} M(x)\nabla \partial^3_{t}u(x,2t)\cdot \overline{\W(x,2t)}\,dx\right| 
\,\leq\, K(t) + L(t)\,,
\]
where, by (2.7) and the definition of $f_t$,
\begin{equation*}
K(t):=
\left|\int_{\rz}M(x)\left(\nabla  \partial^2_{t}\mathcal{D}_t\right)(f_t)(x)\cdot  \overline{\W(x,2t)}\,dx\right|\,
 =\,\left| \int_{\rz}\partial_{t}u(\cdot, t) \,
\overline{\left(\partial_{\nu_{A^*}}\partial_t^2(\mathcal{S}^{\L^*}_{-t}\nabla)\right)
\big(M^*\W(\cdot,2t)\big)}\,dx\right|\,,
\end{equation*}
and by (2.6)
\begin{multline*}
L(t):=
\left|\int_{\rz}M(x) \left(\nabla  \partial^2_{t} \big(\mathcal{S}_t\nabla_{\parallel}\big)\right) 
 \left( \big(A\nabla u(\cdot, t)\big)_{\parallel} \right)(x) \cdot 
\overline{ \W(x,2t)}\,dx\right| \\[4pt]
=\,  \left|\int_{\rz}\big(A\nabla u(\cdot, t)\big)_{\parallel}\cdot \overline{
\left(\nabla_{\parallel}\partial_t^2(\mathcal{S}^{\L^*}_{-t}\nabla)\right)(M^*\W(\cdot,2t))}\,dx\right|.
\end{multline*}
In turn, plugging these bounds into (3.9) and using Cauchy-Schwarz, we obtain the bound
\begin{multline*}
I_1' \lesssim \int_0^\infty \big(K(t) + L(t)\big) \,t dt \\[4pt]
\lesssim
\l(\int_{0}^{\infty}\!\int_{\rz}\big\vert t \nabla 
u(x, t)\big\vert^2\frac{dxdt}{t}\r)^{1/2}\l(\int_{0}^{\infty}\!\int_{\rz}
\bigg\vert t^3 \left(\nabla\partial_t^2\big(\mathcal{S}^{\L^*}_{-t}\nabla\big)\right)\l(M^* \W(\cdot,2t)\r)(x)
\bigg\vert^2\frac{dxdt}{t}\r)^{1/2} \\[4pt]
=:\, {\bf A}\times {\bf B}\,.
\end{multline*} 
Observe that ${\bf A} =  \interleave  t\nabla u \interleave$, and since $u$ is a solution,
by the case $m=1$ of (2.12), we have
\[
{\bf A} = \interleave  t\nabla u \interleave\lesssim \sup_{t>0}\|u(\cdot, t)\|_{L^2(\rz)}\,.
\]
Consider now the factor {\bf B}.  Recall that our goal is to show that
${\bf B}\lesssim \|N_*(\W)\|_{L^2(\rz)} +
\interleave  t\nabla \W \interleave$.
To this end, we set 
\[
\Theta_t:= t^3 \nabla\partial_t^2\big(\mathcal{S}^{\L^*}_{-t}\nabla\big)\,,
\]
and as above, let $P_t$ be a nice approximate identity with a smooth,
radial, compactly supported kernel.
We then write 
\[
\Theta_t\big(M^* \W(\cdot,2t)\big)(x) \,=\, \Theta_t \, M^*(x) \, P_t  \W(\cdot,2t)(x) \, +\, 
\R_t  \W(\cdot,2t)(x) \,,
\]
where 
for a function ${\bf f}$ valued in $\mathbb{C}^{n+1}$ (in particular, for ${\bf f}= W(\cdot,2t)$ with $t$
momentarily fixed), we define
\[
\R_t \,{\bf f} (x) \, := \,\Theta_t\left(M^* {\bf f}\right)(x)\, - \,\Theta_t\,M^*(x)\, P_t \,{\bf f}(x)\,.
\]
(For future reference, we observe that one may define
$\R_t$ on $(n+1)\times(n+1)$ matrix-valued
functions in the obvious way).
We then have 
\[
{\bf B} \, \leq \,    \interleave  
\left(\Theta_t \, M^*\right) \, \big(P_t  \W(\cdot,2t) \big)\interleave \,+\,
\interleave   \R_t  \W(\cdot,2t) \interleave.
\]
By the Carleson measure estimate
(2.10), applied to $\L^*$, in the lower half space, with $m=3$, we have 
\[
\interleave  \left(\Theta_t\, M^* \right)\, \big(P_t \W(\cdot,2t)\big) \interleave
\lesssim \|N_{*}(P_t\W(\cdot,2t))\|_{L^2(\rz)} 
\lesssim\, \|N^3_{*}(\W)\|_{L^2(\rz)} 
\,\lesssim\, \|N_{*}(\W)\|_{L^2(\rz)}\,, 
\]
where in the last step we have used the well-known observation of \cite{FS} that 
non-tangential maximal functions defined using cones with different apertures
are equivalent in $L^p$ norm.

Finally, we consider the contribution of the
remainder term $\R_t$.
By (2.11) applied to $\L^*$, in the lower half-space, with $m=3$,
\[
\|\Theta_t\, (\f1_{2^{k+1}Q\setminus 2^kQ})\|^2_{L^2(Q)}
\lesssim 
2^{-nk}2^{2k}\l(\frac{t}{2^k \ell(Q)}\r)^{6}\|\f\|^2_{L^2(2^{k+1}Q\setminus 2^kQ)}
\lesssim 
2^{-nk}\l(\frac{t}{2^k \ell(Q)}\r)^{4}\|\f\|^2_{L^2(2^{k+1}Q\setminus 2^kQ)}\,, \eqno(3.12)
\]
uniformly for each $k\geq 1$, and $0<t\leq 16\,\ell(Q)$.
In addition, by (2.8), $\Theta_t$ is bounded on $L^2(\rn,\mathbb{C}^{n+1})$, uniformly in $t>0$.
By \cite[Lemma 3.11]{AAA} and the definition of $P_t$,
these facts continue to hold with $\R_t$ in place of $\Theta_t$.
In turn, this 
allows one to define $\R_t {\bf 1}$ as an
element of $L^2_{\loc }$, where ${\bf 1}$ denotes the
$(n+1)\times(n+1)$ identity matrix,
and by construction $\R_t {\bf 1}=0$, since the approximate identity
$P_t$ preserves constants.
Thus, we may apply \cite[Lemma 3.5]{AAA}
to $\R_t$, to deduce that
 $$\int_{\rz}\vert \R_t\W(x,2t)\vert^2dx\lesssim t^2 \int_{\rz}\vert\nabla_x \W(x,2t)\vert^2dx.$$ 
Consequently,
\[
\interleave   \R_t  \W(\cdot,2t) \interleave \lesssim 
\interleave   t\nabla  \W(\cdot,2t) \interleave\,,
\]
as desired.
This completes the proof of Lemma 3.1, and hence that of Theorem 1.1. 
\hfill$\Box$

\section{Solvability with $L^\infty$ data, and an Agmon-Miranda Maximum Principle}

Recall the following result of Auscher \cite{A}:

\begin{theorem}[\cite{A}]
Let $A$, $A_0$ and $\L$ be as above, but possibly $t$-dependent. If $\ez$ is small enough, depending only on $n$ and $\gamma$, then there is a positive exponent $\alpha$ and a constant $C$ (each depending only on $n$ and $\gamma$) such that, given $u$ solving $\L u=0$ in a ball $2B:=B(X, 2R)$, with $R>0$, $$|u(Y)-u(Z)|\leq C\l(\frac{|Y-Z|}{R}\r)^{\alpha}\l(\fint_{2B}|u|^{2}\r)^{1/2}, \quad \forall\; Y, Z \in B=B(Y, R). \eqno(4.1)$$ (Here, capital letters 
denote points in $\rdl$, e.g., $X:=(x,t)$).
\end{theorem}
	
	From  Theorem 4.1, we may deduce the following.

\noindent{{\bf Corollary 4.2.} {\it Let $A$, $A_0$ and $\L$ be as above, but possibly $t$-dependent. If $\ez$ is small enough, depending only on $n$ and $\gamma$, then there is a positive exponent $\alpha$ and a constant $C$ (each depending only on $n$ and $\gamma$) such that, given any 
cube $Q \subset \rz$, and its double $2Q$, 
along with their associated Carleson boxes $R_Q:=Q\times (0, l(Q))$, and 
$R_{2Q}:=2Q\times (0, 2l(Q))$, and a solution $u\in W^{1,2}(R_{2Q})$, vanishing 
in the trace sense on $2Q$, 
then $$|u(x,t)|\leq C\l(\frac{t}{\ell(Q)}\r)^{\alpha}\l(\frac{1}{\ell(Q)^{n+1}}
\iint_{R_{2Q}}|u|^{2}\r)^{1/2}, \quad \forall\; (x,t) \in R_{Q}. \eqno(4.2)$$ }

\noindent
{\it Proof.} The proof 
follows immediately from Theorem 4.1 by making an 
odd reflection across the boundary $2Q\times \{0\}$. We omit the details. 
 \hfill$\Box$

     \medskip                                                                                                                                                                              
                                                                                                                                                                                   
\noindent{{\bf Corollary 4.3.} {\it Let $A$, $A_0$ and $\L$ be as 
in Theorems 1.1, 1.2 (in particular, $t$-independent), and
4.1, with $\ez$ small enough that $(D)^{\L}_2$ and $(R)^{\L}_2$ are both solvable (see Theorem 2.2). 
Let $f \in L^2(\rz)$, and let $u$ be the solution of $(D)^{\L}_2$ with boundary data $f$. 
If $f$ vanishes on $2Q$, then the 
conclusion of Corollary 4.2 continues to hold. }

\medskip

\noindent
{\it Proof.} Note that if we were to
assume $f\in H^1(\rn)$, then 
the solution of $(R)^{\L}_2$ with boundary data $f$ 
satisfies the assumptions in Corollary 4.2, thus (4.2) holds.
Moreover, as  we have previously mentioned,
by \cite{HK}, the problems $(D)^{\L}_2$ and 
$(R)^{\L}_2$ are compatibly solvable, in particular, for data $f\in H^1(\rz)$, 
the solution of $(D)^{\L}_2$ with data $f$ equals the solution of 
$(R)^{\L}_2$ with data $f$ (the latter is unique only up to an additive constant, 
but will be equal to the former for a suitable choice of this constant).

Since $f \in L^2(\rz)$, and vanishes on $2Q$, we can 
approximate $f$ in $L^2$ norm by $f_k\in C^{\infty}\cap H^1(\rn)$,
with each $f_k$ vanishing on $\frac{3}{2}Q$.   
Let $u_k$ denote the solution to $(D)^\L_2$, and compatibly, to $(R)^\L_2$,
with data $f_k$. 
Since Corollary 4.2 clearly holds with $\frac{3}{2}Q$ in place of  $2Q$, 
we find that (4.2) holds for each $u_k$, uniformly in $k$.
We may then pass to the limit as follows.
Observe that (4.2) holds with $u$ 
replaced by $u_k$, and that by the $L^2$ estimates for $(D)^{\L}_2$, 
$$\sup_{t>0}\|u(\cdot , t)-u_k(\cdot , t)\|_2\lesssim
 \|f-f_k\|_2\rightarrow 0, \quad \mbox{as}\;  k \rightarrow \infty.$$ 
Consequently, for $(x,t)\in \rdd$ fixed, combining the latter estimate with
the interior Moser-type local boundedness estimate we obtain
\begin{multline*}
|u(x, t)-u_k(x, t)|\, \lesssim\, \l(\fint_{t/2}^{2t}\fint_{|x-y|<t}|u(y, s)-u_k(y, s)|^2\, dyds\r)^{1/2} \\[4pt]
 \lesssim \,  \l(t^{-n}\fint_{t/2}^{2t}\int_{\rn}|u(y, s)-u_k(y, s)|^2\, dyds\r)^{1/2} \, \to\, 0\,,
 \quad \mbox{as}\;  k \rightarrow \infty.
\end{multline*}
Similarly, for any fixed cube $Q\subset \rn$,
\[
\iint_{R_{2Q}}|u-u_k|^{2} \leq \int_0^{\ell(Q)}\!\!\int_{\rn} |u(y,s)-u_k(y,s)|^{2}\, dy ds 
 \, \to\, 0\,,
 \quad \mbox{as}\;  k \rightarrow \infty.
\]
 We conclude that (4.2) holds for $u$.
\hfill$\Box$

In the sequel, let $$\Delta(x, r):=\{y\in \rz: |x-y|<r\}$$ denote 
the surface ball of radius $r$ and center $x$, on $\rn \cong \partial \rdd$.

\smallskip

\noindent{{\bf Lemma 4.4.} {\it Let $A,\L$ be as in Corollary 4.3.
Let $x\in \rz$, and $0<t<\frac{1}{100}R$, with $R\leq R'<\infty.$ 
Suppose that $g\in L^{\infty}$ with 
$$\supp(g)\subset S_{R, R'}=S_{R, R'}(x):=\Delta(x, R')\setminus \Delta(x, R).$$ 
Let $v$ solve $(D)^{\L}_2$ with boundary data $g$. 
Then there exists a constant $C=C(n, \gamma)$ such that $$|v(x, t)|\leq C\l(\frac{t}{R}\r)^{\alpha}\|g\|_{L^{\infty}},\eqno(4.3)$$ uniformly in $R'$, for $R'\geq R$, where $\alpha>0$ is the exponent in Corollaries 4.2 and 4.3. }

\medskip

\noindent
{\it Proof.} Set $$\Delta_k=\Delta_k(x):=\Delta(x, 2^k):=\{y\in \rz: |x-y|<2^k\}, \;k=0, 1, 2, ...,$$ and 
$$ 
S_k=\Delta_{k+1}\setminus\Delta_k, \quad k\in \mathbb{Z}.$$ Thus 
$$S_{R, R'}\subset \bigcup_{k=k(R)}^{k(R')}S_k,$$ where 
$2^{K(R)}\approx R$ and $2^{K(R')}\approx R'$. 
Consequently, $$v=\dsum_{k=k(R)}^{k(R')}v_k,$$ 
where $v_k$ solves $(D)^{\L}_2$ with boundary data $g_k:=g1_{S_k}$. 
By Corollary 4.3 and the solvability of $(D)^{\L}_2$, we have 
$$|v_k(x, t)|\lesssim \l(\frac{t}{2^k}\r)^{\alpha}\l(\fint_{0}^{2^k}\!\fint_{\Delta_k}|v_k|^{2}\r)^{1/2}
\lesssim\l(\frac{t}{2^k}\r)^{\alpha}\l(2^{-kn}\int_{S_k}|g|^2\r)^{1/2}\lesssim\l(\frac{t}{2^k}\r)^{\alpha}\|g\|_{L^{\infty}}.$$ 
Summing up $k\geq k(R)$, we get (4.3).

 \hfill$\Box$
 
We are now able to prove the following.

\noindent{{\bf Proposition 4.5.} {\it Let $A$, $A_0$ and $\L$ be as in Corollary 4.3.
Let $f\in L^{\infty}(\rz)$. Then there is a solution $u$ of $\L u = 0$ in $\rdd$ such that $u(\cdot, 0) = f$ in the sense of non-tangential convergence, satisfying the Agmon-Miranda maximum principle $$\|u\|_{L^{\infty}(\rdd)}\leq C\|f\|_{L^{\infty}(\rz)},$$ where $C=C(n, \gamma)$. }

\medskip

\noindent
{\it Proof.}\quad Given a point $x_0\in \rz$, we define the dyadic surface balls centered at $x_0$ on 
$\rz \cong \rz \times \{0\}=\partial \rdd$ by 
$$\Delta_k=\Delta_k(x_0):=\Delta(x_0, 2^k):=\{y\in \rz: |x_0-y|<2^k\}, \;k\in \zc$$ 
and set 
\[
S_k=\Delta_{k+1}\setminus\Delta_k, \quad k\in \mathbb{Z}\,,
\]
so that $\cup_kS_k=\rz\setminus \{0\}$.

We let $f_k:=f1_{S_k}$ and $u_k$ be the solution of $(D)^{\L}_2$ with boundary data $f_k$. 
Define $$f^{N}:=\dsum_{k=-\infty}^{N}f_k, \quad u^N:=\dsum_{k=-\infty}^{N}u_k.$$ 
Clearly, $u^N$ is the unique solution
of $(D)^{\L}_2$ with boundary data $f^N$. 
To prove the proposition, we will show that $$u:=\lim_{N\rightarrow \infty}u^N$$ exists at 
each point of $\rdd$ and satisfies the conclusion of  the theorem. Moreover, $u$ is well-defined, 
in the sense that if $u'$ is 
constructed in the same way as $u$, but for a different center $x'_0$, then $u = u'.$ To this end, 
we fix a point $(x, t)\in \rdd$ and suppose that 
$2^M\geq 2^N\gg t+|x-x_0|$. Then by the definition of $f^N$, 
$$\supp(f^M-f^N)\subset S_{R, R'}(x)\,,$$ 
where $R\approx 2^N$ and $R'\approx 2^M$. 
By Lemma 4.4, 
$$|u^N(x, t)-u^M(x, t)|\leq C\l(\frac{t}{2^N}\r)^{\alpha}\|f\|_{L^{\infty}(\rz)}\rightarrow 0,
\quad \mbox{as}\; N, M \rightarrow 0.\eqno(4.4)$$
Thus, $u^N$ converges pointwise, and in fact, uniformly on compacta, in $\rdd$, hence also in $L^2_{loc}(\rdd)$. 
By Caccioppoli's inequality applied to 
$u^N-u^M$, we further 
see that $u^N$ converges in $W^{1,2}_{loc}(\rdd)$, 
whence the limit $u$ also solves $\L u=0$ in $\rdd.$

Let us now show that $u$ satisfies the required properties.

\smallskip

\noindent
{\it Definition of $u$ is independent of center $x_0$.} Fix two distinct points $x_1, x_2 \in \rz$ and construct the corresponding $f^{N, i}, u^{N, i}, i=1, 2$, and $u^i=\lim_{N\rightarrow \infty}u^{N, i}, i=1,2$ as above, with 
$x_1, x_2$ in place of $x_0$. Let $(x, t)\in \rdd$ and consider $M, N$ such that
 $$2^M\geq 2^N\gg t+|x-x_1|+x-x_2|\,.$$ 
 Then $$\supp(f^{M, 1}-f^{N, 2})\subset S_{R, R'},$$ where 
 $R\approx 2^N$ and $R'\approx 2^M$. Again we invoke Lemma 4.4 to get 
 $$|u^{N, 2}(x, t)-u^{M, 1}(x, t)|\leq C\l(\frac{t}{2^N}\r)^{\alpha}\|f\|_{L^{\infty}(\rz)}\rightarrow 0,
 \quad \mbox{as}\; N, M \rightarrow 0.\eqno(4.5)$$ Therefore, $u^{N, 2}$ and $u^{M, 1}$ converge to the 
 same limit $u$.

\smallskip

\noindent
{\it Non-tangential convergence to $f$.} Fix $x_0\in \rn$, and build $f^N$ and $u^N$ as above, 
relative to the center $x_0$.
By $(D)^{\L}_2$, each $u^N$ converges 
non-tangentially to $f^N$, for a.e. $x \in \rz$. Thus, there is a set $Z  =\cup_N Z_N\subset \rz$, of 
measure zero, such that  
$u^N$ converges non-tangentially to $f^N$ for every $N$, and
at every point $x \in \rz\setminus Z$.
 Fix such an $x\neq x_0$ and let $\ez>0$. Consider the truncated cone at $x$ of 
 height $\ez:$ 
 $$\Gamma^{\ez}(x):=\{(y, t)\in \rdd: |x-y|<t<\ez\}.$$ 
 Observe that (4.4) 
 continues to hold with $y$ in place of $x$, for $(y, t) \in \Gamma^{\ez}(x)$, with $\ez$ small, and $N, M$ large. We therefore have for such $(y, t)$ that $$|u(y, t)-u^N(y, t)|=\lim_{M\rightarrow \infty}|u^M(y, t)-u^N(y, t)|\lesssim \l(\frac{\ez}{2^N}\r)^{\alpha}\|f\|_{L^{\infty}(\rz)}.\eqno(4.6)$$ On the other hand, if we fix $N$ so large that 
 $f^N(x)=f(x)$ and use that $u^N$ converges non-tangentially to $f$ at $x$, then for 
 $(y, t)\in \Gamma^{\ez}(x)$, we have 
 $$|u^N(y, t)-f(x)|=o(1), \quad \mbox{as}\; \ez\rightarrow 0.$$ 
 Letting $\ez\rightarrow 0$, we see that $u(y, t)\rightarrow f(x)$ non-tangentially.

\smallskip

\noindent
{\it Agmon-Miranda maximum principle.}\quad Let $(x, t)\in \rdd$. 
We seek 
to show that
 $$|u(x, t)|\leq C\|f\|_{L^{\infty}}\,,$$ 
 with $C=C(n, \gamma)$. Since the definition of $u$ is
  independent of the choice of $x_0$ used in the 
  construction, we may choose $x_0 = x.$ 
  We then define 
  $\Delta_k=\Delta_k(x),\, S_k=S_k(x), \,f_k=f1_{S_k}$ 
  and $u_k$ 
  as above. Choose $k(t)$ such that $2^{k(t)}\approx t$ and write 
  $$u=\dsum_{k=-\infty}^{k(t)+10}u_k+\dsum_{k=k(t)+11}^{\infty}u_k=:U_1+U_2.$$ 
  By Moser local boundedness and
$(D)^{\L}_2$, 
$$|U_1(x,t)|\lesssim\l(\fint_{t/2}^{2t}\fint_{|x-y|<t}|U_1(y,s)|^2\r)^{1/2}\lesssim 
\l(t^{-n}\int_{|x-y|<C2^{k(t)}}|f|^2\r)^{1/2}\lesssim \|f\|_{L^{\infty}},$$ 
since $2^{k(t)}\approx t$, where all of the implicit constants in the display depend only on dimension and ellipticity. Furthermore, for $k>k(t)$, by Lemma 4.4, we have $$|u_k(x, t)|\lesssim \l(\frac{t}{2^k}\r)^{\alpha}\|f\|_{L^{\infty}(\rz)},$$ and so we may sum over $k\geq k(t)+11$ to see that $|U_2(x,t)|\lesssim \|f\|_{L^{\infty}(\rz)}.$
\hfill$\Box$

\smallskip

\noindent{\bf Remark 4.}
Note that by construction, if $f \in L^\infty(\rn)$ is {\em compactly supported}, then the solution of
$(D)^{\L}_{\infty}$ with boundary data $f$, and the solution of
$(D)^{\L}_{2}$ with boundary data $f$, are the same.

\smallskip

We conclude this section with the following.  Recall that $\Delta(x,r)$ denotes the ``surface ball" 
centered at $x$, of radius $r$, on
$\rn\cong \partial \rdd$.  Given $\Delta=\Delta(x,r)$, let $R_\Delta:= \Delta\times (0,r)\subset \rdd$
denote the usual Carleson cylinder above $\Delta$.

\smallskip

\noindent{{\bf Proposition 4.6.} {\it Let $A$, $A_0$ and $\L$ be as in Proposition 4.5 (i.e., as in
Corollary 4.3).
Let $f\in L^{\infty}(\rz)$, and let $u \in L^\infty(\rdd)$ be the solution 
of $\L u=0$ in $\rdd$, with data $f$, constructed 
in Proposition 4.5.   Set $d\mu(x,t):= |t\nabla u(x,t)|^2\,t^{-1} dxdt$.
We then have the Carleson measure estimate
\[
\|\mu\|_{\mathit{c}}\,:=\,\sup_\Delta\frac{1}{|\Delta|} 
\iint_{R_\Delta} |t\nabla u(x,t)|^2\, \frac{dxdt}{t} \, \leq \, C \|f\|^2_{L^\infty(\rn)}\,,
\]
where $C$ depends only on dimension and ellipticity.
}

\noindent
{\it Proof.} Given our preceding work in this section, 
the argument is standard, but we include it here for the sake of completeness.  Fix a 
surface ball $\Delta_0:=\Delta(x_0,r)\subset \rn$, set $\Delta_k:= \Delta(x_0,2^k r)$,
and $S_k:=\Delta_{k+1}\setminus \Delta_k$.  Now define $f_k=f1_{S_k}$, let $u_k$ solve
$(D)^{\L}_{\infty}$ (equivalently, $(D)^{\L}_{2}$, see Remark 4) 
with boundary data $f_k$, and as in the proof of Proposition 4.5, set
$$f^{N}:=\sum_{k=-\infty}^{N}f_k, \quad u^N:=\sum_{k=-\infty}^{N}u_k,
\eqno(4.7)$$ 
so that $u^N$ is the solution
of $(D)^{\L}_2$ (and of $(D)^{\L}_{\infty}$) with boundary data $f^N$. 
As noted above, $u^N \to u$ in $W^{1,2}_{loc}(\rdd)$, hence, for each $\delta\in(0,r)$,
\[
\int_\delta^r \!\!\int_{\Delta_0} |t\nabla u^N(x,t)|^2\, \frac{dxdt}{t} \,\to\,
\int_\delta^r \!\!\int_{\Delta_0} |t\nabla u(x,t)|^2\, \frac{dxdt}{t} \,,\quad \text{as } N\to\infty.
\]
Thus, it is enough to show that
\[
r^{-n} \int_0^r \!\!\int_{\Delta_0} |t\nabla u^N(x,t)|^2\, \frac{dxdt}{t} \,\lesssim \, 
\|f\|^2_{L^\infty(\rn)}\,, \eqno(4.8)
\]
uniformly in $N$.  Using the notation of (4.7), we write
\[
u^N \,=\, \sum_{k=-\infty}^{0}u_k \,+ \,\sum_{k=1}^{N}u_k \,=\, u^0 \,+ \,\sum_{k=1}^{N}u_k\,,
\]
so that
\begin{multline*}
\left(\int_0^r \!\!\int_{\Delta_0} |t\nabla u^N(x,t)|^2\, \frac{dxdt}{t}\right)^{1/2}\\[4pt]
\leq\,
\left(\int_0^r \!\!\int_{\Delta_0} |t\nabla u^0(x,t)|^2\, \frac{dxdt}{t}\right)^{1/2}
\,+\, \sum_{k=1}^{N}
\left(\int_0^r \!\!\int_{\Delta_0} |t\nabla u_k(x,t)|^2\, \frac{dxdt}{t}\right)^{1/2}
\,=:\, I_0 \,+\, \sum_{k=1}^{N} I_k\,.
\end{multline*}
By (2.12) with $m=1$, and the solvability of $(D)^{\L}_{2}$, we have
\[
(I_0)^2 \, \lesssim \, \|f^0\|^2_{L^2(\rn)} \,=\, 
\int_{\Delta(x_0,2r)} |f(x)|^2\, dx\,
\lesssim\, r^{n}  \|f\|^2_{L^\infty(\rn)}\,,
\]
as desired.  

By construction, $f_k$ vanishes outside of $S_k$, so
by Corollary 4.3, $u_k$ is H\"older continuous up to the boundary outside of $S_k$, 
and we may therefore
use Caccioppoli's inequality at the boundary and then Corollary 4.3 
(i.e., inequality (4.2), but with surface balls in place of cubes), to write 
\begin{multline*}
(I_k)^2 \lesssim r \int_0^r \!\!\int_{\Delta_0} |\nabla u_k(x,t)|^2\, dxdt\,
\lesssim\, r^{-1} \int_0^{2r} \!\!\int_{\Delta(x_0,2r)} | u_k(x,t)|^2\, dxdt\\[4pt]
\lesssim \, 2^{-2k\alpha} r^n\|u_k\|^2_{L^\infty(\rdd)}\,
\lesssim \, 2^{-2k\alpha} r^n \|f\|^2_{L^\infty(\rn)}\,,
\end{multline*}
where in the last step we have used the Agmon-Miranda maximum principle.
We may now sum a geometric series to conclude.
\hfill$\Box$

\section{Proof of Theorem 1.2}

In this section, we focus on the proof of Theorem 1.2, which, together with the 
results in the previous section, comprise
 the main new contributions of this paper. 
The proof will be split into two parts. In Part 1, 
we present a suitable definition 
of the commutator $[\Lambda, g](f)$, 
under the assumptions that $f\in L^{\infty}(\rz)$ and 
$g\in \dot{H}^1(\rz)$.   In Part 2, we prove a variant of
 Dahlberg's bilinear estimate by a more refined version of the procedure used to prove 
 Lemma 3.1. The conclusion of the theorem then follows.
As in the preceding section, we let
$\Delta(x, r):=\{y\in \rz: |x-y|<r\}$ denote 
the surface ball of radius $r$ and center $x$, on $\rn$.

\smallskip

\noindent{{\bf Part 1: definition of  $\|[\Lambda, g](f)\|_{L^2}$.}}

Under the hypotheses of Theorem 1.2, we have from Theorem 2.2 that both $(D)_2$ and $(R)_2$ are solvable for $\L$,  and its adjoint $\L^*$.

We let $\Lambda^*$ to denote the adjoint of  
$\Lambda$.  Observe that $\Lambda^*$ is the the Dirichlet-to-Neumann map for the adjoint operator $\L^*$,
as may be seen by the Gauss-Green formula. 

For $f\in L^{\infty}(\rz)$, let $u$ be the solution of $(D)^{\L}_{\infty}$ with boundary data $f$, 
as constructed in section 4.  
We may assume that $g\in C_0^{\infty}(\rz)$, which is dense in $\dot{H}^1(\rz)$. 

For $0<\delta\ll 1$, and $1\ll R<\infty$, set $f_{\delta}:=u(\cdot, \delta)$, and choose $\eta_{R}\in C_0^{\infty}(B(0, 2R))$ with $\eta_{R}\equiv 1$ on $B(0, R)$. Let $u_{\delta, R}$ be the solution of $(D)^{\L}_2$
(equivalently, the solution of $(D)^{\L}_\infty$; see Remark 4) 
with boundary data $f_{\delta}\eta_{R}$. By the Agmon-Miranda maximum principle proved 
in section 4, $f_{\delta}:=u(\cdot, \delta)$ satisfies that 
$$\lim_{\delta\rightarrow 0}f_{\delta}(x)=f(x), \;\mbox{a.e.}\; 
x\in \rz\,,\quad \mbox{and}\quad\sup_{\delta>0}\|f_{\delta}\|_{L^{\infty}(\rz)}\lesssim 
\|f\|_{L^{\infty}(\rz)}.\eqno(5.1)$$
For any $h\in C_0^{\infty}(\rz)$, we shall prove the following estimate 
$$|I_{\delta,R}|:=\biggl|\int_{\rz}\l[\Lambda(f_{\delta}\eta_{R}g)\overline{h}-g\Lambda(f_{\delta}\eta_{R})
\overline{h}\r]\biggl|\lesssim \|f_{\delta}\|_{L^{\fz}}\|h\|_{L^2}\|\nabla_\| g\|_{L^2}\lesssim
 \|f\|_{L^{\fz}}\|h\|_{L^2}\|\nabla_\| g\|_{L^2},\eqno(5.2)$$ 
 uniformly in $\delta$ and $R$; in fact, the implicit constants depend only on $n, \gamma$, provided $\ez$ is small enough, with the same dependence. Observe that we have used (5.1) in the last step.

Taking (5.2) for granted momentarily, 
we seek to extend estimate (5.2) to the limiting case as 
$\delta\rightarrow 0$ and $R\rightarrow \infty$. 
To this end, we define $$\int_{\rz}\l[\Lambda(fg)\overline{h}-g\Lambda(f)\overline{h}\r]:=\lim_{ R\rightarrow \infty}\lim_{\delta\rightarrow 0}\int_{\rz}\l[\Lambda(f_{\delta}
\eta_{R}g)\overline{h}-g\Lambda(f_{\delta}\eta_{R})\overline{h}\r]
=:\lim_{ R\rightarrow \infty}\lim_{\delta\rightarrow 0}I_{\delta,R}.\eqno(5.3)$$
Let us show that this definition is reasonable, and in particular that the limit exists. We observe that at least formally,
$$\int_{\rz}\l[\Lambda(fg)\overline{h}-g\Lambda(f)\overline{h}\r]=\int_{\rz}\l[fg\overline{\Lambda^*(h)}-f\overline{\Lambda^*(\overline{g}h)}\r],\eqno(5.4)$$ so our goal is to show that 
the limit in (5.3) exists, and is equal to the right hand of (5.4).

By  \cite[Theorem 1.14]{AAA} (solvability of $(R)_2^\L$), 
the analogue of (5.4) does hold for any $\delta>0$ and $R<\infty$, i.e., we can write $I_{\delta,R}$ as $$I_{\delta,R}=\int_{\rz}\l[f_{\delta}\eta_{R}\l(g\overline{\Lambda^*(h)}-\overline{\Lambda^*(\overline{g}h)}\r)\r].\eqno(5.5)$$ By the solvability of  $(R)^{\L}_2$, we know that 
$\Lambda^*(h), \Lambda^*(\overline{g}h)\in L^2$ (recall that we have taken $g,h\in C^\infty_0$ by density).
Consequently, by (5.1) we may use dominated convergence to obtain
$$\lim_{\delta\rightarrow 0}I_{\delta,R}=\lim_{\delta\rightarrow 0}\int_{\rz}\l[f_{\delta}\eta_{R}\l(g\overline{\Lambda^*(h)}-\overline{\Lambda^*(\overline{g}h)}\r)\r]=\int_{\rz}\l[f\eta_{R}\l(g\overline{\Lambda^*(h)}-\overline{\Lambda^*(\overline{g}h)}\r)\r]:=I_R.$$ Since (as we shall prove) (5.2) holds uniformly in 
$\delta>0$, we also have that 
$$|I_R|\lesssim \|f\|_{L^\infty}\|h\|_{L^2}\|\nabla_\| g\|_{L^2}.\eqno(5.6)$$ 
Set $$\Psi=\Psi(g, h):=\l(g\overline{\Lambda^*(h)}-\overline{\Lambda^*(\overline{g}h)}\r).$$ Since (5.6) holds for any $f\in L^{\infty}(\rz)$, we have 
$$\sup_{R<\infty} \int_{|x|<R} |\Psi(x)| \,dx
\,\lesssim \,\|h\|_{L^2}\|\nabla_\| g\|_{L^2},$$ 
by the definition of $\eta_{R}$, and thus using monotone convergence theorem we also have that $$\|\Psi\|_{L^1(\rz)}\lesssim \|h\|_{L^2}\|\nabla g\|_{L^2}.$$ Consequently, we obtain the desired limit 
$$\lim_{ R\rightarrow \infty}\lim_{\delta\rightarrow 0}I_{\delta,R}=
\lim_{ R\rightarrow \infty}\int_{\rz}f\eta_R\Psi(g, h)=\int_{\rz}f\Psi(g, h)$$ 
by dominated convergence theorem. This completes Part 1.
It remains to prove (5.2).

\medskip


\noindent{{\bf Part 2: the proof of (5.2).}}

We now fix $0<\delta\ll 1$, and $1\ll R<\infty$, let $f_\delta$ and $\eta_R$ be defined as in Part 1 above,
and for notational convenience,
we set $f = f_\delta \eta_R$. 
Recall also that by density, we may assume that $g,h \in C_0^\infty$.
Then qualitatively, with this revised notation, $fg\in H^1(\rz)$ and $f\in H^1(\rz)$. 
Of course, by hypothesis, we also have a quantitative $L^\infty$ bound for $f$, and moreover
$f$ now has compact support.
We let $u$ be the solution of
$(D)^{\L}_{\infty}$ (equivalently, the solution of
$(D)^{\L}_{2}$; see Remark 4) with boundary data $f$, and as above, we let
$H$ be the solution of $(D)^{\L^*}_2$ with boundary data $h$,
and set  $V(x, t)=\varphi_t*g =P_t\,g$, where $P_t$ is a nice approximate identity
with a smooth, radial, compactly supported kernel $\varphi_t$. 

Thus, (5.2) will follow immediately 
once we establish the following estimate:
$$\biggl|\int_{\rz}\l[\Lambda(fg)\overline{h}-g\Lambda(f)\overline{h}\r]\biggl|
\,\lesssim\, \|f\|_{L^{\fz}(\rn)}\,\|h\|_{L^2(\rn)}\,\|\nabla_\| g\|_{L^2(\rn)}\,,$$
for all $g,h \in C_0^\infty(\rn)$, and every 
 $f \in H^1(\rn)\cap L^\infty(\rn)$ with
 compact support.

Exactly as in $(3.3)$,  we have  
$$\int_{\rz}\l[\Lambda(fg)\overline{h}-g\Lambda(f)\overline{h}\r]
\,=\,\iint_{\rdd}u\nabla V\overline{A^*\nabla H}\,-\,\iint_{\rdd}A\nabla u \nabla V\overline{H}:=I+J\,.
\eqno(5.7)$$
By Lemma 3.1,  and the solvability of 
$(D)^{\L^*}_{2}$, 
\begin{multline*}
|I|\, \lesssim \, \sup_{t>0}\|H(\cdot,t)\|_{L^2(\rn)}\,
\bigg(\|N_{*}(u\nabla V)\|_{L^2}
\,+\,\interleave  t\nabla (u \nabla V) \interleave\bigg)\\[4pt]
\lesssim\,
\|h\|_{L^2}\l(\|N_{*}(u\nabla V)\|_{L^2}\,+\,
\interleave  t\nabla u\cdot \nabla V \interleave\,+\,\interleave  tu \nabla^2 V \interleave\r)\\[4pt]
=:\, \|h\|_{L^2}\,\big({\bf M}_1 + {\bf M}_2+{\bf M}_3\big)\,.
\end{multline*}
In turn, to handle term $I$, it is therefore enough to show that
\[
{\bf M}_1 + {\bf M}_2+{\bf M}_3 \, \lesssim \,  \|f\|_{L^{\fz}(\rn)}\, \|\nabla_\| g\|_{L^2}\,.
\]
To this end, let us note that
by Proposition 4.6,
$$d\mu=\vert t\nabla u(x,t)\vert^2\frac{dxdt}{t} \;\mbox{is a Carleson measure on }\;\rdd \;\mbox{with norm}\; 
\|\mu\|_{\mathit{c}}\lesssim \|f\|^2_{L^{\fz}(\rn)}\,.\eqno(5.8)$$ 
Recall that $V=P_t g$, so that by (3.2), $N_*(\nabla V) 
\lesssim \mathcal{M} (\nabla_\| g)$, where
$\mathcal{M}$ denotes the Hardy-Littlewood maximal operator.   Moreover, 
$t\nabla^2 V = Q_t(\nabla_\| g)$, where $Q_t$ satisfies the classical Littlewood-Paley estimate
\[
\interleave \! Q_t F \interleave \lesssim \|F\|_{L^2(\rn)},
\]
for arbitrary $F\in L^2(\rn)$.
Consequently,
\[
\|N_{*}(\nabla V)\|_{L^2}\, +\, 
\interleave  t \nabla^2 V \interleave\,\lesssim  \,\|\nabla_{\parallel} g\|_{L^2} \eqno(5.9)
\]
With these observations in hand, by the Agmon-Miranda maximum principle and (5.9), we have
\[
{\bf M}_1\,\leq\, \|u\|_{L^{\fz}(\rdd)}\,\|N_{*}(\nabla V)\|_{L^2}\, 
 \lesssim \,  \|f\|_{L^{\fz}(\rn)}\, \|\nabla_\| g\|_{L^2}\,.
\]
By (5.8), Carleson's lemma, and (5.9),
\[
{\bf M}_2 \,\lesssim\, \|f\|_{L^{\fz}(\rn)}\,\|N_{*}(\nabla V)\|_{L^2}\, 
 \lesssim \,  \|f\|_{L^{\fz}(\rn)}\, \|\nabla_\| g\|_{L^2}\,,
\] 
and by the Agmon-Miranda maximum principle and (5.9), 
\[
{\bf M}_3  \, \leq \,
\|u\|_{L^{\fz}(\rdd)}\interleave  t \nabla^2 V \interleave\, 
\lesssim\,
 \|f\|_{L^{\fz}(\rn)}\, \|\nabla_\| g\|_{L^2}\,.
\]
This concludes out treatment of term $I$.

It remains to 
estimate term $J$ (see (5.7) above), which is the heart of the matter.
The basic strategy will be that of Lemma 3.1, but in the present setting we shall
need to proceed more carefully. 
 As in the proof of Lemma 3.1, it suffices to prove 
$$\dsup_{0<\rho\ll 1}\fint_{\rho}^{2\rho}\bigg\vert\int^{1/\theta}_{\theta}\!
\!\int_{\rz}A(x)\nabla u(x,t)\cdot \nabla V(x, t) \,\overline{H(x,t)}\, dxdt\bigg\vert d\theta
\,\lesssim \,\|f\|_{L^{\fz}(\rn)}\, \|h\|_{L^2(\rn)} \,\|\nabla_{\parallel} g\|_{L^2(\rn)}.$$
For any $\rho>0$ small, integrating by parts in $t$, we have the following 
\begin{multline*}
\dint^{1/\theta}_{\theta}\int_{\rz}A(x)\nabla u(x,t)\cdot \nabla V(x, t) 
\overline{H}dxdt \\[4pt]
=\,\left(\dint_{\rz} A\,t\nabla u(\cdot,t)\cdot \nabla V(\cdot, t)
 \overline{H}\,dx\right)\bigg|_{t=\theta}^{t=1/\theta}
\,-\, \dint^{1/\theta}_{\theta}\int_{\rz}  A\nabla \partial_{t}u\cdot \nabla V\,
\overline{H}\,dx\, tdt \\[4pt]
- \, \dint^{1/\theta}_{\theta}\int_{\rz}  A\, t\nabla u\cdot \nabla \partial_{t}V\,
\overline{H}\,dx\,dt\ 
\,-\, \dint^{1/\theta}_{\theta}\int_{\rz}  A\, t\nabla u \cdot \nabla V\,
\overline{ \partial_{t}H}\,dx\, dt  \\[4pt]
\,=:\, J_{1}-J_{2}-J_{3}-J_{4}.
\end{multline*} 

We start with the last of these.  Uniformly in $\theta$, and hence in $\rho$, we have
\begin{multline*} 
|J_{4}|\,\lesssim\, \l(\dint_{\rdd}\big|t \nabla u(x,t)\cdot \nabla V(x, t)
\big|^2\frac{dxdt}{t}\r)^{1/2}\l(\int_0^\infty \!\int_{\rn} 
\big|t\partial_t H(x,t)\big|^2\, \frac{dxdt}{t}\r)^{1/2}\\[4pt]
\lesssim \, 
\interleave  t\nabla u\cdot \nabla V \interleave \interleave t \nabla H \interleave 
\,\lesssim \,\|f\|_{L^{\fz}(\rn)}\|\nabla_{\parallel} g\|_{L^2}  \|h\|_{L^2},
\end{multline*} 
where we have used (5.8), (5.9) (as for the term ${\bf M}_2$ above), 
along with (2.12) for the adjoint solution $H$, and  the solvability of $(D)^{\L^*}_{2}$.

By Cauchy-Schwarz,
\begin{multline*}
|J_{3}|\,\lesssim \, \interleave \, |t \nabla u| \,\,|H(\cdot,t)| \,\interleave
\interleave  t \nabla^2 V \interleave
\\[4pt]
\lesssim\,  \|f\|_{L^{\fz}(\rn)}\,\|N_*H\|_{L^2(\rn)}\,\|\nabla_\| g\|_{L^2(\rn)}\,\lesssim\,
 \|f\|_{L^{\fz}(\rn)}\,\|h\|_{L^2(\rn)}\,\|\nabla_\| g\|_{L^2(\rn)}\,,
\end{multline*} 
as desired, uniformly in $\theta$ (hence also in $\rho$),
where we have used (5.8) and Carleson's lemma, (5.9), and the fact that $H$ solves
$(D)_2^\L$ with data $h$. 

The boundary terms $J_{1}$ are handled as follows:
\begin{multline*}
\fint_{\rho}^{2\rho}|J_{1}|\,d\theta \,\lesssim\,
\sup_{r>0} \fint_{r}^{2r} \!\!\int_{\rz} |t \nabla u(x,t) |\, |\nabla V(x, t)| \,
|H(x,t)|\,dx \,dt\\[4pt]
\lesssim\, \left(\int_0^\infty  \!\!\int_{\rz} |t \nabla u(x,t) |^2\, |\nabla V(x, t)|^2\,\frac{dxdt}{t}  \right)^{1/2}
\|N_*(H)\|_{L^2(\rn)}\\[4pt]
\lesssim\,
 \|f\|_{L^{\fz}(\rn)}\,\|\nabla_\| g\|_{L^2(\rn)}\,\|h\|_{L^2(\rn)}\,,
\end{multline*} 
uniformly in $\rho$,
by (5.8), Carleson's lemma and (5.9), and the fact that $H$ solves
$(D)_2^{\L^*}$ with data $h$.

It remains to treat $J_{2}$.  To this end, we begin by recording the following 
generalization of
(5.8), which follows from the latter by the $t$-independence of $A$ and Caccioppoli's inequality in 
Whitney boxes:
for any $m\geq 1$, 
$$d\mu_m=\vert t^{m}\nabla \partial_{t}^{m-1} u(x,t)\vert^2\frac{dxdt}{t} 
\;\mbox{is a Carleson measure on }\;\rdd \;\mbox{with norm}\; 
\|\mu_m\|_{\mathit{c}}\lesssim \|f\|^2_{L^{\fz}(\rn)},\eqno(5.10)$$ 
with implicit constant depending of course on $m$, as well as on dimension and ellipticity.

To control the term $J_{2}$,  
we integrate by parts up to a total of $N+1$ times in $t$ (that is, $N$ additional times:  
we have already done so once), 
for some suitably large integer $N$ to be chosen,
stopping the first time that a $t$-derivative falls on 
either $\nabla V$ or $\overline{H}$.
In either of the latter two cases, the result is a term of
the same form as $J_{3}$ or $J_{4}$, along with boundary terms of the same form as $J_{1}$,
except with $t\nabla u$ replaced by $t^{m} \nabla \partial_t^{m-1} u$, for some $2\leq m\leq N+1$.
Using (5.10) in lieu of (5.8), we may handle these terms exactly like their counterparts with $m=1$,
already treated above.  The one scenario that remains to be considered is that which occurs when all
$N+1$ derivatives in $t$ fall upon $u$, i.e., it remains only to show that 
\begin{multline}\tag{5.11}
\dsup_{0<\rho\ll 1}\fint_{\rho}^{2\rho}\bigg|\dint^{1/\theta}_{\theta}\!\!\int_{\rz} A(x)
\nabla \partial_{t}^{N+1}u(x,t)\cdot \nabla V(x, t) \overline{H(x,t)}\,dx\,  t^{N+1} dt\bigg|d\theta \\[4pt]
=: \, \dsup_{0<\rho\ll 1}\fint_{\rho}^{2\rho}|\Omega(\theta)| \,d\theta \, \lesssim \,
\|f\|_{L^{\fz}(\rn)} \|h\|_{L^2(\rn)} \|\nabla_{\parallel} g\|_{L^2(\rn)} , 
\end{multline}
provided that $N$ is chosen large enough;
in particular, it will be enough to take $N=n+2$ in the sequel. 

To prove (5.11), we shall follow the outline of the 
argument in Section 3.  We first make the change of variable $t\rightarrow 2t$,
to obtain
\[
\Omega(\theta)\,=\, C_N \dint^{1/\theta}_{\theta}\!\!\int_{\rz} A(x)
\nabla \partial_{t}^{N+1}u(x,2t)\cdot \nabla V(x, 2t) \overline{H(x,2t)}\,dx\,  t^{N+1} dt\,,
\]
and then we use the Green formula (3.10) (bearing in mind our qualitative assumptions on $u$),
and set $s=t$ , to get the following generalization of (3.11): 
\[
\left(\nabla \partial^{N+1}_{t}u\right)(\cdot,2t)=-\left(\nabla  \partial^N_{t}\mathcal{D}_t\right)(f_t)\,-\,
\left(\nabla  \partial^N_{t} \big(\mathcal{S}_t\nabla_{\parallel}\big)\right) 
\left(\big(A\nabla u(\cdot, t)\big)_{\parallel}\right)\,, 
\]
where as before, $f_t := \partial_t u(\cdot,t)$.
We may then set $\W(\cdot,t):= \overline{\nabla V(\cdot,t)} H(\cdot,t)$, and  
use (2.6) and (2.7) to write
\[
\left|\int_{\rz} A(x)
\nabla \partial_{t}^{N+1}u(x,2t)\cdot \nabla V(x, 2t) \overline{H(x,2t)}\,dx\right|\,
\leq \, K_N(t) + L_N(t)\,,
\]
where
\begin{equation*}
K_N(t):=
\left|\int_{\rz}A(x)\left(\nabla  \partial^N_{t}\mathcal{D}_t\right)(f_t)(x)\cdot  \overline{\W(x,2t)}\,dx\right|\,
 =\,\left| \int_{\rz}\partial_{t}u(\cdot, t) \,
\overline{\left(\partial_{\nu_{A^*}}\partial_t^N(\mathcal{S}^{\L^*}_{-t}\nabla)\right)
\big(A^*\W(\cdot,2t)\big)}\,dx\right|\,,
\end{equation*}
and 
\begin{multline*}
L_N(t):=
\left|\int_{\rz}A(x) \left(\nabla  \partial^N_{t} \big(\mathcal{S}_t\nabla_{\parallel}\big)\right) 
 \left( \big(A\nabla u(\cdot, t)\big)_{\parallel} \right)(x) \cdot 
\overline{ \W(x,2t)}\,dx\right| \\[4pt]
=\,  \left|\int_{\rz}\big(A\nabla u(\cdot, t)\big)_{\parallel}\cdot \overline{
\left(\nabla_{\parallel}\partial_t^N(\mathcal{S}^{\L^*}_{-t}\nabla)\right)(A^*\W(\cdot,2t))}\,dx\right|.
\end{multline*}
Using (2.9), we observe that these expressions make sense, by virtue of our qualitative assumptions on
$u$, and the fact that $H(\cdot,t)\in L^\infty(\rn)$
 (qualitatively, because the data $h\in C_0^\infty$; see Remark 4), for each fixed $t>0$,
and therefore $\W(\cdot,t) \in L^2(\rn)$ (again qualitatively).
Note that
\[
K_N(t) + L_N(t) \,\lesssim\, \int_{\rn} \big|\nabla u(\cdot, t)\big| \, 
\left| \left(\nabla\partial_t^N(\mathcal{S}^{\L^*}_{-t}\nabla)\right)(A^*\W(\cdot,2t))\right|\,dx\,,
\]
hence, plugging this bound into the definition of $\Omega(\theta)$, and in turn
into  (5.11), it suffices to prove that
\[
\int_0^\infty\!\!\int_{\rn} \big|t\nabla u(\cdot, t)\big| \, 
\left| \Theta_t\left(A^*\W(\cdot,2t)\right)\right|\,\frac{dxdt}{t}\,
\lesssim \, \|f\|_{L^{\fz}(\rn)}\, \|h\|_{L^2(\rn)}\, \|\nabla_{\parallel} g\|_{L^2(\rn)} \,, \eqno(5.12)
\]
where
\[
\Theta_t:=  t^{N+1}\nabla \partial_{t}^N (\mathcal{S}^{\L^*}_{-t}\nabla)
\]
(note that the operator $\Theta_t$ defined in Section 3 was exactly the same, but with $N=2$).
Let $P_t$ be the nice  approximation of the identity 
with a smooth, compactly supported kernel, introduced previously.
Just as in Section 3, we then write 
\[
\Theta_t\big(A^* \W(\cdot,2t)\big)(x) \,=\, \Theta_t \, A^*(x) \, P_t  \W(\cdot,2t)(x) \, +\, 
\R_t  \W(\cdot,2t)(x) \,,
\]
where 
for a function ${\bf f}$ valued in $\mathbb{C}^{n+1}$ (in particular, for ${\bf f}= W(\cdot,2t)$ with $t$
momentarily fixed), we define
\[
\R_t \,{\bf f} (x) \, := \,\Theta_t\left(A^* {\bf f}\right)(x)\, - \,\Theta_t\,A^*(x)\, P_t \,{\bf f}(x)\,.
\]

We first consider the contribution of $ \Theta_t \, A^*(x) \, P_t  \W(\cdot,2t)(x)$ in (5.12).
Note that  
\begin{multline*}
 \frac{1}{|Q|}\int_0^{l(Q)}\!\!\int_{Q}\big| t \nabla u(x,t)\big| \, \big|\Theta_t(A^*)\big| \, 
\frac{dxdt}{t}\\[4pt]
\leq \,  \l(\frac{1}{|Q|}\int_0^{l(Q)\!\!}\int_{Q}\big|t \nabla u(x,t)\big|^2\frac{dxdt}{t}\r)^{1/2}
\l(\frac{1}{|Q|}\int_0^{l(Q)}\!\!\int_{Q}\big|\Theta_t(A^*)\big|^2\frac{dxdt}{t}\r)^{1/2}
 \\[4pt]
 \lesssim\,  \|f\|_{L^\infty(\rn)}
  \, \|A\|_{L^\infty(\rn)}\, \approx \,  \|f\|_{L^\infty(\rn)}\,,
\end{multline*}
uniformly in $Q$, by (2.10), (5.8) and ellipticity. 
Recall that $\W(\cdot,t):= \overline{\nabla V(\cdot,t)} H(\cdot,t)$, so
 by Carleson's lemma, we have
\begin{multline*}
\int_0^\infty\!\!\int_{\rn} \big|t\nabla u(x, t)\big| \, 
\left|  \Theta_t \, A^*(x) \, P_t  \W(\cdot,2t)(x) \right|\,\frac{dxdt}{t}\,\\[4pt]
\lesssim \, \|f\|_{L^{\fz}(\rn)}\, \|N_*\big(P_{t}\W(\cdot, 2t)\big)\|_{L^1(\rn)} 
\lesssim \, \|f\|_{L^{\fz}(\rn)}\,
\|N_*(\nabla V)\|_{L^2}\|N_*H\|_{L^2} \\[4pt]
\lesssim \, \|f\|_{L^{\fz}(\rn)}\, \|h\|_{L^2(\rn)}\, \|\nabla_{\parallel} g\|_{L^2(\rn)}\,,
\end{multline*}
as desired.

Last, we deal with the remainder term $\R_t$. We begin by recording two facts, for future reference. 
The first 
entails precise quantitative dependence on the aperture of the cones used to define
the non-tangential maximal function:
$$\|N_{*}^{\beta}(f)\|_{L^2}\, \lesssim \beta^{n/2}\|N_{*}(f)\|_{L^2}\eqno(5.13)$$
for any $f\in L^2$ and $\beta\geq 1$; the proof can be found in \cite[Lemma 1, p. 166]{FS}. 
The second indicates the off-diagonal decay for $\Theta_t$, and hence for $\R_t$: 
for every cube $Q$ and all $t\lesssim \ell(Q)$, 
$$\|\R_t(\f1_{2^{j+1}Q\setminus 2^jQ})\|^2_{L^2(Q)}\lesssim 
2^{-nj}2^{2j}\l(\frac{t}{2^j \ell(Q)}\r)^{2N+2}\|\f\|^2_{L^2(2^{j+1}Q\setminus 2^jQ)}, 
\quad \forall j\geq 1,\eqno(5.14)$$
for any $\f\in L^2(\rz,\ccl)$. For $\Theta_t$, the latter estimate is simply (2.11) for $\L^*$, in the lower
half-space, with $m=N+1$. 
As in Section 3, where we considered the case $N=2$ (see (3.12) above), we may use
(2.8), \cite[Lemma 3.11]{AAA} and the definition of $P_t$, to 
extend the estimate to $\R_t$, which is (5.14).
As in Section 3, we may then define $\R_t {\bf 1}$ as an
element of $L^2_{\loc }$, where ${\bf 1}$ denotes the
$(n+1)\times(n+1)$ identity matrix,
and by construction $\R_t {\bf 1}=0$.

As above, let $\mathbb{D}_k$ denote the grid of dyadic cubes on $\rn$ of length $\ell(Q) = 2^k$.
Let $Q\in \mathbb{D}_k$, suppose that $t\in(2^{k}, 2^{k+1}]$, and for $i\geq 1$, set
$[\W]_{2^iQ}:= \fint_{2^iQ} \W(\cdot,2t)$.
For $j\geq 1$, since $t\approx \ell(Q)$,
we then have 
\begin{multline*}
\left(\dint_{2^{j+1}Q} \big|\W(x, 2t)-[\W]_{2Q}\big|^2dx\right)^{1/2} \\[4pt]
\lesssim  \left(\dint_{2^{j+1}Q}\big|W(x, 2t)-[\W]_{2^{j+1}Q}\big|^2dx\right)^{1/2}
\,+\, 
\dsum_{i=1}^j
\left(2^{jn}|Q|\, \big|[\W]_{2^{i+1}Q} - [\W]_{2^iQ}\big|^2\right)^{1/2}\\[4pt]
\lesssim\, 2^{j} \left(\int_{2^{j+1}Q}
|t\nabla_{x}\W(x, 2t)|^2dx\right)^{1/2}\,
+\, \sum_{i=1}^j\left(2^{(j-i)n}\, 2^{2i}\dint_{2^{i+1}Q}|t\nabla_{x}W(x, 2t)|^2dx\right)^{1/2}
\\[4pt]
\lesssim\, j\,2^{jn/2} 2^{j}\, \left(\int_{2^{j+1}Q}
|t\nabla_{x}\W(x, 2t)|^2dx\right)^{1/2}\,,
\end{multline*}
by Poincar\'e's inequality.
Thus, since $\R_t{\bf 1} = 0$, and $t\approx \ell(Q)$, we see from (5.14) that 
\begin{multline}\tag{5.15}
\left(\int_{Q} 
\left|  \R_t  \W(\cdot,2t)(x) \right|^2\, dx\right)^{1/2} \, 
\leq \, 
\left( \int_{Q}  \, 
\left|  \R_t \, \bigg( \big(\W(\cdot,2t)(x) - [\W]_{2Q}\big) 1_{2Q}\bigg)\right|\, dx\right)^{1/2} \\[4pt]
\qquad\qquad\qquad \qquad\qquad \qquad
+\, \,\,
\sum_{j=1}^\infty
\left( \int_{Q}  \, 
\left|  \R_t \, \bigg( \big(\W(\cdot,2t)(x) - [\W]_{2Q}\big) 1_{2^{j+1}Q\setminus 2^jQ}\bigg)\right|\, dx\right)^{1/2}
\\[4pt]
\lesssim\,  \sum_{j=1}^\infty j\,2^{-j(N-1)} \left(\int_{2^{j+1}Q}
|t\nabla_{x}\W(x, 2t)|^2dx\right)^{1/2}
\end{multline}
We shall now use the preceding estimate to establish the following.

\smallskip

\noindent{\bf Claim}.  Define the conical square function
\[
\A \W(x):= \left(\iint_{|x-y|<t} \big|\R_t\W\big(\cdot,2t\big)(y)\big|^2\, \frac{dydt}{t^{n+1}}\right)^{1/2}.
\]
We then have
\[
\|\A\W\|_{L^1(\rn)} \, \lesssim \, \|\nabla_\|g\|_{L^2(\rn)}\, \| h\|_{L^2(\rn)}\,. \eqno(5.16)
\]

\noindent{\it Proof of Claim}.  Using (5.15), we find that for some purely dimensional constant $M$,
\begin{multline*}
\A\W(x)\, \leq \,
\left(\sum_{k=-\infty}^\infty \, \sum_{Q\in \mathbb{D}_k:\, \text{dist}(x,Q)<2^{k+1}}\,
\int_{2^k}^{2^{k+1}}\!\!\!\int_{Q} \big|\R_t\W\big(\cdot,2t\big)(y)\big|^2\, \frac{dydt}{t^{n+1}}\right)^{1/2}
\\[4pt]
\lesssim \, \sum_{j=1}^\infty j\,2^{-j(N-1)}
\left(\sum_{k=-\infty}^\infty \, \sum_{Q\in \mathbb{D}_k:\, \text{dist}(x,Q)<2^{k+1}}\,
\int_{2^k}^{2^{k+1}}\!\!\!\int_{2^{j+1}Q}
|t\nabla_{y}\W(y, 2t)|^2\, \frac{dydt}{t^{n+1}}\right)^{1/2}\\[4pt]
\lesssim\, 
\sum_{j=1}^\infty j\,2^{-j(N-1)}
\left(
\iint_{|x-y|< M 2^j t}
|t\nabla_{y}\W(y, 2t)|^2\, \frac{dydt}{t^{n+1}}\right)^{1/2}\,.
\end{multline*}
Recall that in the present context,
$\W(\cdot,t) = \overline{\nabla V(\cdot,t)} H(\cdot,t)$.
For notational convenience, we set 
\[
g_1:= \nabla_\| g, \quad G_1(x,t):=  \nabla V(\cdot,2t)(x),  \qquad 
g_2:= h,\quad G_2(x,t):=  H(x,2t)\,,
\]
so that
\[ |t\nabla_\|\W(\cdot,2t)|\,\lesssim\, |t\nabla G_1|\, |G_2| \, +\, |t\nabla G_2|\, |G_1|\,.
\]
Note that by (5.9), (2.12), and the solvability of $(D)_2^{\L^*}$,
\[
\|N_{*}(G_i)\|_{L^2(\rn)}\, +\, 
\interleave  t \nabla G_i \interleave\,\lesssim  \,\| g_i\|_{L^2(\rn)}\, \qquad i=1,2\,.
\]
Thus, to prove the claim, it suffices to show that
\[
 \sum_{j=1}^\infty j\,2^{-j(N-1)} \int_{\rn}
\left(
\iint_{|x-y|< M 2^j t}
|t\nabla G_1(y,t)|^2\,|G_2(y,t)|^2 \frac{dydt}{t^{n+1}}\right)^{1/2} dx
\,\lesssim \interleave  t \nabla G_1 \interleave\, \|N_{*}(G_2)\|_{L^2(\rn)} \,, \eqno(5.17)
\]
along with a similar estimate with the roles of $G_1$ and $G_2$ reversed.
Since the roles of $G_1$ and $G_2$ are symmetrical, we need only treat the version
stated in (5.17).  Note that for $|x-y|<M2^jt$,
 we have
 \[|G_2(y,t)| \leq N_*^{M2^j} G_2(x)\,,
 \]
 i.e., the non-tangential maximal function defined with respect to a cone of aperture
 $M2^j$.  Thus, the left hand side of (5.17) is bounded by
 \begin{multline*}
  \sum_{j=1}^\infty j\,2^{-j(N-1)} \int_{\rn} N_*^{M2^j} G_2(x)\,
\left(
\iint_{|x-y|< M 2^j t} 
|t\nabla G_1(y,t)|^2\, \frac{dydt}{t^{n+1}}\right)^{1/2} dx\\[4pt]
\leq \,
 \sum_{j=1}^\infty j\,2^{-j(N-1)} \,  \|N_*^{M2^j} G_2\|_{L^2(\rn)}\,
\left( \int_{\rn} 
\iint_{|x-y|< M 2^j t} 
|t\nabla G_1(y,t)|^2\, \frac{dydt}{t^{n+1}} dx\right)^{1/2}\\[4pt]
\lesssim \, 
 \sum_{j=1}^\infty j\,2^{-j(N-1)} \,2^{jn/2} \, \|N_* G_2\|_{L^2(\rn)}\, 2^{jn/2}\, 
\interleave  t \nabla G_1 \interleave\,,
 \end{multline*}
 where in the last step we have used (5.13), along with the following 
estimate, obtained via Fubini's theorem:
 \[
  \int_{\rn} 
\iint_{|x-y|< M 2^j t} 
|t\nabla G_1(y,t)|^2\, \frac{dydt}{t^{n+1}} dx\,
=\,
 \int_0^\infty\!\!\int_{\rn} |t\nabla G_1(y,t)|^2\, 
t^{-n}\int_{|x-y|< M 2^j t} 
dx\, \frac{dydt}{t}\, \approx\, 2^{jn} \interleave  t \nabla G_1 \interleave^2.
 \]
 We now choose $N=n+2$, to obtain (5.17), and hence the claim.
\hfill$\Box$

With (5.16) in hand, and using the Carleson measure estimate (5.8),
we then obtain
\[
 \|\mu\|^{1/2}_{\mathit{c}}\,\|\A\W\|_{L^1(\rn)}
\lesssim \|f\|_{L^{\fz}(\rn)} \|\nabla_\|g\|_{L^2(\rn)}\, \| h\|_{L^2(\rn)}\,.\eqno(5.18)
\]
We also claim that
\[
{\bf K}\, :=\,
\int_0^\infty\!\!
\int_{\rn} \big|t\nabla u(y, t)\big| \, 
\left|  \R_t  \W(\cdot,2t)(y) \right|\, \frac{dydt}{t}\,\lesssim \, \|\mu\|^{1/2}_{\mathit{c}}\,\|\A\W\|_{L^1(\rn)}\,.
\eqno(5.19)
\]
Momentarily taking (5.19) for granted, we 
then immediately obtain the desired estimate (5.12) for the contribution of the $\R_t$ term, 
by combining (5.18)-(5.19).  The conclusion of Theorem 1.2 follows.

It remains only to discuss (5.19).  In fact, the latter is actually
a classical estimate of Fefferman
(see \cite[pp. 148-149]{FS}), but for the reader's convenience, we shall reproduce
the argument here.  To this end, for $0<h< \infty$, set
\[
\A_h \W(x):= \left(\iint_{|x-y|<t\leq h} \big|\R_t\W\big(\cdot,2t\big)(y)\big|^2\, 
\frac{dydt}{t^{n+1}}\right)^{1/2},\quad
\A_h(t\nabla u) (x):= \left(\iint_{|x-y|<t\leq h} \big|t\nabla u(y,t)\big|^2\, \frac{dydt}{t^{n+1}}\right)^{1/2}\,.
\]
(thus, for all $h\in (0,\infty)$, $\A_h\W \leq \A\W$ as defined above).
By (5.8) (i.e., Proposition 4.6), for all $y\in \rn$, and all $h\in (0,\infty)$,
\[
\int_{|y-x|<h} \big(\A_h(t\nabla u) (x)\big)^2 dx \, \leq \, C_0  \|\mu\|_{\mathit{c}} \, h^n\,,\eqno(5.20)
\]
with $C_0$ depending only on dimension.
Set
\[
h(x):=\sup\left\{h\geq 0: \, \A_h(t\nabla u) (x) \leq C_1  \|\mu\|^{1/2}_{\mathit{c}} \right\}\,,
\]
with $C_1$ a sufficiently large dimensional constant to be chosen momentarily.
Note that in particular,
\[
\A_{h(x)}(t\nabla u) (x) \leq C_1 \|\mu\|^{1/2}_{\mathit{c}}\,.\eqno(5.21)
\]
Then for every $y\in \rn$,  there is a uniform constant $c$ such that
\[
|\{x\in \rn:\, |x-y|<h \leq h(x)\}| \geq c h^n\,. \eqno(5.22)
\]
Indeed, by definition, if $h(x)<h$, then $\A_h(t\nabla u) (x) >C_1  \|\mu\|^{1/2}_{\mathit{c}}$,
so that by Tchebychev's inequality
\begin{multline*}
\left|\big\{x:\, |x-y|<h \text{ and }  h> h(x)\big\}\right| \leq 
\left|\big\{x:\, |x-y|<h \text{ and } \A_h(t\nabla u) (x) >C_1  \|\mu\|^{1/2}_{\mathit{c}}\big\}\right|\\[4pt]
\leq \frac{1}{C_1^2 \|\mu\|_{\mathit{c}}} \int_{|x-y|<h} \left( \A_h(t\nabla u) (x)\right)^2\, dx \,
\leq \, \frac12 \left|\big\{x:\, |x-y|<h \big\}\right|\,,
\end{multline*}
by (5.20), provided that $C_1$ is chosen large enough, depending on $C_0$.
Consequently, using (5.22), we see that
\begin{multline*}
{\bf K} \, \lesssim \,
 \int_0^\infty\!\!
\int_{\rn} \big|t\nabla u(y, t)\big| \, 
\left|  \R_t  \W(\cdot,2t)(y) \right|\, t^{-n} \int_{|x-y|<t<h(x)} dx\, \frac{dydt}{t}\\[4pt]
= \int_{\rn} \iint_{|x-y|<t<h(x)}  \big|t\nabla u(y, t)\big| \, 
\left|  \R_t  \W(\cdot,2t)(y) \right|\,\frac{dydt}{t^{n+1}}\, dx\\[4pt]
\lesssim \,  \int_{\rn} \left(\iint_{|x-y|<t<h(x)}  \big|t\nabla u(y, t)\big|^2 
\,\frac{dydt}{t^{n+1}}\right)^{1/2}\,
\left(\iint_{|x-y|<t<h(x)} 
\left|  \R_t  \W(\cdot,2t)(y) \right|^2\,\frac{dydt}{t^{n+1}}\right)^{1/2} dx\\[4pt]
\leq\,  \int_{\rn} \A_{h(x)}(t\nabla u) (x) \, \A \W(x) \, dx\, \lesssim\, 
 \|\mu\|^{1/2}_{\mathit{c}} \,\|\A\W\|_{L^1(\rn)}\,,
\end{multline*}
by (5.21), so that (5.19) holds.

This concludes the proof of Theorem 1.2.
\hfill$\Box$


\medskip

S. Hofmann,\quad
  
Department of Mathematics, University of Missouri, Columbia, MO 65211, USA

E-mail address:  hofmanns@missouri.edu

G. Zhang,\quad

School of Mathematical Sciences, Peking University, Beijing, 100871, P. R. China

E-mail address:  zhangguoming256@pku.edu.cn
\end{document}